\documentclass[11pt]{article}

\usepackage{cite}
\usepackage{amsmath,amssymb,amsfonts}
\usepackage{amsthm,xparse}
\usepackage{algorithmic}
\usepackage{graphicx}
\usepackage{textcomp}
\usepackage{xcolor}
\usepackage{soul}
\usepackage{url}
\usepackage{bm}

\makeatletter
\DeclareRobustCommand{\thmprime}{%
  \begingroup
  \expandafter\in@\expandafter b\expandafter{\f@series}%
  \ifin@ \boldmath \fi
  $\m@th{}^{\prime}$%
  \endgroup
}
\DeclareRobustCommand{\thmprimeprime}{%
  \begingroup
  \expandafter\in@\expandafter b\expandafter{\f@series}%
  \ifin@ \boldmath \fi
  $\m@th{}^{\prime\prime}$%
  \endgroup
}

\newcommand*\rel@kern[1]{\kern#1\dimexpr\macc@kerna}
\newcommand*\widebar[1]{%
  \begingroup
  \def\mathaccent##1##2{%
    \rel@kern{1.2}% 0.8 both
    \overline{\rel@kern{-1.25}\macc@nucleus\rel@kern{-0.1}}% 0.2 both
    \rel@kern{0.1}%
  }%
  \macc@depth\@ne
  \let\math@bgroup\@empty \let\math@egroup\macc@set@skewchar
  \mathsurround\z@ \frozen@everymath{\mathgroup\macc@group\relax}%
  \macc@set@skewchar\relax
  \let\mathaccentV\macc@nested@a
  \macc@nested@a\relax111{#1}%
  \endgroup
}
\makeatother

\newtheorem{theorem}{Theorem}
\newtheorem{proposition}[theorem]{Proposition}
\newtheorem{lemma}[theorem]{Lemma}
\newtheorem{corollary}[theorem]{Corollary}
\theoremstyle{remark}

\theoremstyle{definition}

\newtheorem{problem}[theorem]{Problem}

\NewDocumentEnvironment{manual}{O{theorem}m}
 {
  \addtocounter{theorem}{-1}%
  \renewcommand\thetheorem{#2}%
  \begin{#1}
 }
 {\end{#1}}

\def\secure{secure}

\def\uniform{uniform }

\def\cdc{,\ldots,}
\def\N{{\mathbb N}}
\def\Nodd{{\N_{\mathrm{odd}}}}

\def\calG{{\mathcal{G}}}
\def\calD{{\mathcal{D}}}
\def\calGnd{\calG_{n|d}}
\def\calGndh{\calG_{n|d|h}}
\def\hmin{h_{\min}}
\def\hmax{h^-_{\max}}
\def\hmind{\hmin(n,d)}
\def\hminn{\hmin(n)}
\def\hmand{\hmax(n,d)}

\def\bH{{\mathbf H}}
\def\bS{{\mathbf S}}

\def\inline#1:{\par\vskip 7pt\noindent{\bf #1:}\hskip 10pt}

\def\lI{\!\!\!\!\!}
\def\lJ{\!\!\!\!\!\!\!}
\def\lK{\!\!\!\!\!\!\!\!\!\!}
\def\lLl{\!\!\!\!\!\!\!\!\!\!\!\!\!\!\!\!\!\!\!\!\!\!\!\!\!\!\!\!\!\!\!\!\!\!\!\!\!}
\def\lL{\!\!\!\!\!\!\!\!\!\!\!\!\!\!\!\!\!\!\!\!\!\!\!\!\!\!\!\!\!\!\!\!\!\!\!\!\!\!\!\!}
\def\maj{\gamma_{\mathrm{maj}}}
\def\smaj{\gamma_{\mathrm{smaj}}}
\def\q{k}
\def\qq{k}
\def\kk{q}

\def\sta#1{\left(#1\right)^*}
\def\staa#1{\left(#1\right)^{**}}
\newcommand{\x}[1]{} %\x{} For in-line comments

\def\xy{\hspace{.07em}}
\newcommand{\eq}[2]{\begin{equation}\label{#1}#2\end{equation}}
\newcommand{\eqss}[1]{\begin{eqnarray}#1\end{eqnarray}}
\newcommand{\conne}[2]{\big[#1\!\equiv\!\equiv\!#2\big]}
\newcommand{\connee}[3]{\big[#1\!\equiv\!\equiv\!#2\!\equiv\!\equiv\!#3\big]}
\newcommand{\doubling}[2]{\big[#1{\,\succ\hspace{-5pt}-\,}#2\big]}   %>\hspace{-5pt}- \nolinebreak
\newcommand{\mediastinum}[3]{\big[#1{\,\succ\hspace{-5pt}-\,}#2\!\equiv\!\equiv\!#3\big]}  %>
\newcommand{\composiTwo}[2]{\big[S_{#1}\!\equiv\!\equiv\!H_{#2}\big]}
\newcommand{\composiThree}[3]{\big[S_{#1}\!\equiv\!\equiv\!S_{#2}\!\equiv\!\equiv\!H_{#3}\big]}

\newcommand{\card}[1]{|#1|}
\newcommand{\lU}[1]{#1^{\mathstrut}}
\newcommand{\lD}[1]{#1_{\mathstrut}}
\newcommand{\lB}[1]{#1^{\mathstrut}_{\mathstrut}}
\newcommand{\indi}[2]{{\bm1}_{\{#1\}}(#2)}
\newcommand{\modd}[2]{(#1\!\!\mod #2)}
\newcommand{\modD}[2]{(#1\!\!\!\!\mod #2)}
\definecolor{ggreen}{rgb}{0,.47,0}
\newcommand{\cha}[2]{#2}

\title{The Power of Small Coalitions\\under Two-Tier Majority on Regular Graphs}
\author{Pavel Chebotarev
\thanks{Technion---Israel Institute of Technology, Haifa 3200003, Israel;
A.A. Kharkevich Institute for Information Transmission Problems of RAS, Moscow 127051, Russia.
Email: {\tt pavel4e@technion.ac.il.}}
\and
David Peleg
\thanks{Department of Computer Science and Applied Mathematics, The Weizmann Institute of Science, Rehovot, Israel. Email: {\tt david.peleg@weizmann.ac.il}.}
}

\begin{document}
	
\maketitle

%\centerline{\large\bf Draft 2.3}
{\abstract
\begin{center}
\begin{minipage}[c]{33em}
In this paper, we study the following problem.
Consider a setting where a proposal is offered to the vertices of a given network $G$, and the vertices must conduct a vote and decide whether to accept the proposal or reject it. Each vertex $v$ has its own valuation of the proposal; we say that $v$ is ``happy'' if its valuation is positive (i.e., it expects to gain from adopting the proposal) and ``sad'' if its valuation is negative. However, vertices do not base their vote merely on their own valuation. Rather, a vertex $v$
is a \emph{proponent\/} of the proposal if \cha{a}{the} majority of its neighbors are happy with it
and an \emph{opponent\/} in the opposite case. At the end of the vote, the network collectively accepts the proposal whenever \cha{a}{the} majority of its vertices are proponents.
We study this problem \cha{on}{for} regular graphs with loops. Specifically, we consider the class $\calGndh$ of $d$-regular graphs of odd order $n$ with all $n$ loops and $h$ happy vertices.
We are interested in establishing
necessary and sufficient conditions for the class $\calGndh$ to contain a labeled graph accepting the proposal, as well as
conditions to contain a graph rejecting the proposal.
We also discuss connections to the existing literature, including that on majority domination, and investigate\x{ some} the properties of the obtained conditions.
\end{minipage}
\end{center}}
\bigskip

%----------------------------------------------------
\section{Introduction}
%----------------------------------------------------

\inline Background and motivation:
A potentially desirable property of decision-making policies, hereafter termed \cha{}{the} \emph{faithfulness}, is that every approved decision is desired by \cha{a}{the} majority of the participating population. More explicitly, consider a proposed decision, whose approval will make a subset $H$ of the population $V$ ``happy'' and the complementary set $S=V\setminus H$ ``sad.''
Then a decision-making policy is \emph{faithful\/} if it ensures that a proposal is approved only if $|H|>|V|/2$.

A straightforward faithful decision-making policy is a purely democratic policy in which individual behavior is selfish. Under this policy, each agent in the population $V$ votes in favor of the proposal if and only if it is happy with it, the votes have equal weight, and the final decision is based on majority. This policy is sometimes referred to as the
{``selfish" or ``egoistic" democratic policy.}

Equally faithful is a democratic voting policy where each agent is aware of the preferences of all other agents, and votes according to the wishes of the majority of the population
(again with votes having equal weight and decision by majority). In this case, the voting will in fact end up with a consensus. This policy is sometimes referred to as the {``altruistic" democratic policy.}

In the current paper we are concerned with two-stage decision-making policies collectively known as ``majority of majorities," where the first stage involves local majority votes carried out inside a number of\x{ small} voting bodies (subsets of the population), and the second stage involves a majority vote among the outcomes obtained at these voting bodies.
In some cases, such a policy can alternatively be viewed as a single-stage voting policy falling ``in-between" the egoistic and altruistic policies, where each agent is aware of the preferences of a local subset of the population (viewed as its ``neighborhood"), and votes according to the wishes of the neighborhood majority.

It has long been known that decisions made by a majority of majorities can actually be supported only by a minority, namely, they are not faithful.\x{Nevertheless,}
Despite this property,\footnote{As mentioned in \cite{Kalai10Noise}, a hierarchical method of voting\x{ formed the basis (underlay)} was at the heart of Lenin’s Utopian concept of centralized democracy implemented in the Soviet Union (the\x{ very} term Soviet itself means council) and its satellites for party institutions, national bodies, and trade unions. ``For party institutions (which were the most important) there could be as many as seven levels.’’} called the \emph{referendum paradox\/} \cite{Nurmi99Paradoxes}, two-stage procedures of this kind continue to be widely used. Political examples include the electors' system of electing the president of the United States and many parliamentary procedures. Political scientists continue to warn that ``if a bare majority of members of a bare majority of (equal-sized) parliaments voted ‘yes’ to a measure, that would mean only slightly more than one quarter of all members voted ‘yes’---and thus a measure could pass with nearly three quarters of members opposed''~\cite{Cooper13bicameral}. Is this effect enhanced or weakened when there are many local voting bodies, they have the same size and can be\x{ strongly} arbitrarily mixed? To answer this question, we\x{ need to} study the majority of majorities on regular graphs.

To formalize this question, let us say that a
procedure for approving
proposals is \emph{$\alpha_1$-protecting}
and \emph{$\alpha_2$-trusting\/} (with $0\le\alpha_1\le\alpha_2\le1$) if it
rejects all proposals for which $\card{H}/\card{V}\le\alpha_1$ and \x{rejects no proposal}accepts all proposals for which $\card{H}/\card{V}>\alpha_2.$
In the intermediate cases where
$\alpha_1<\card{H}/\card{V}\le\alpha_2$, such a procedure may accept or reject the proposal taking into account other factors. In particular, this can be determined by the graph expressing the connections between happy and sad agents.
The classical simple majority procedure is obviously both $1/2$-protecting and $1/2$-trusting.
Expressed in these terms, the question we explore is to determine, {for $d$-regular graphs of order $|V|=n$,} the greatest $\alpha_1$ and the smallest $\alpha_2$ such that the ``majority of majorities" policy is $\alpha_1$-protecting and $\alpha_2$-trusting.

\inline Contributions:
The main result of the paper, Theorem~\ref{t:Classes}, concerns a triple $(n,d,h),$ where $n$ is an odd order of a regular graph with loops, $d$ is the vertex degree, and $h$ is the number of happy vertices (supporting the proposal). The theorem establishes conditions that, for a given triple $(n,d,h)$, determine whether there is a graph whose vertices are labeled as happy or sad with parameters $(n,d,h)$ on which the majority of local (neighborhood) majorities accepts the proposal as well as whether there is such a configuration on which this is not the case.
In addition, we study the properties of the relationship between the final two-stage majority decision and the parameters $n,$ $d,$ and~$h.$

The rest of the paper is organized as follows.
Section~\ref{s:Problem} presents the basic notation and formulation of the problem, and
%Section~\ref{s:Ineq}
provides a simple necessary condition for the acceptance of a proposal by a two-stage majority on a regular graph with all loops.
In Section~\ref{s:Constru},
the main technical statements are proved. These enable to obtain necessary and sufficient conditions for the class of graphs with parameters $n$ and $d$ to contain a graph that accepts (rejects)\x{ by the two-stage majority} a proposal making $h$ vertices happy.
These conditions are gathered in Theorem~\ref{t:Classes} presented in Subsection~\ref{s:Classes}.
In Sections~\ref{s:UncertainGap} and~\ref{s:MaxGap} some properties of the found dependencies are reported. In particular, explicit formulas are obtained for the exact integer values of $d$ that allow the lowest support of accepted proposals and simultaneously allow the highest support for rejected proposals. Section \ref{s:MajorityRegular} presents an alternative formulation of the problem as \emph{majority domination\/} on regular graphs.
Finally, Section~\ref{s:Discus} provides a summary and concluding remarks, as well as a discussion of connections to the previous literature and to
% (various) social
some applications.

%-------------------------------------------------------------------------------------
\section{Preliminaries}
\label{s:Problem}
%----------------------------------------------------
\subsection{The \lowercase{Problem}}

Let $G=(V,E)$ be a graph with vertex set $V=V(G)$ and edge set~$E=E(G);$ $\card{V}=n$. For any $u\in V,$ $N_u=\{v\in V\,|\,(u,v)\in E\}$ is the set of neighbors of~$u.$
We consider graphs $G$ with no multiple edges, but where every vertex $v\in V$ has one loop\x{\footnote{Sometimes, graphs with loops allowed are called \emph{pseudographs.}}} $(v,v)\in E$ (hereafter referred to as \emph{graphs with loops}). Hence\footnote{An equivalent formalism is to consider \emph{closed neighborhoods\/} $N_u\cup\{u\}$ for graphs without loops~\cite{Broere95majority}.} $v\in N_v$ for every $v\in V.$ The \emph{degree\/} of vertex $v$ is $d_v=\card{N_v}.$ This implies that a loop is counted\x{ in the degree} once\footnote{For a discussion of this convention see~\cite[Subsection~5.3]{BoldiVigna02}.}.

We consider a setting where the vertices of a given network $G$ must conduct a vote and decide collectively on whether to accept or reject a\x{ global} proposal that is offered to them.
{Each vertex $v$ has its own preference: it is \emph{happy} (respectively, \emph{sad}) with the proposal if it expects to gain (resp., lose) from accepting it. The preference is expressed by a private \emph{opinion} (or \emph{valuation$)$ function} $f:V\to\{1,-1\}$ such that $f(v)=1$ if and only if $v$ is happy with the proposal.}
Thus, an opinion function $f$ induces a \emph{configuration\/} in the form of a binary-vertex-labeled graph $G^f=(V,E,f)$.
\x{Such graphs $G^f$ will be called \emph{binary labeled} graphs.}
The configuration can alternatively be described in terms of a partition of $V$ induced by $f$, dividing it into two disjoint parts $V=\bH\cup \bS$ ($\bH\cap \bS=\varnothing$), where $\bH$ and $\bS$ are called the set of \emph{happy} vertices and the set of \emph{sad} vertices, respectively. Throughout, we denote $|\bH|=h$ and $|\bS|=s$.

However, the choice made by each vertex $v$ during the voting process is not based merely on its own valuation $f(v)$. Rather, it relies on the valuations of $v$'s neighbors as well.
For any $W\subseteq V$, let $f(W)=\sum_{v\in W}f(v)$. In particular, $f(V)$ is called the \emph{weight\/} of~$f$.
We say that a vertex $v\in V$ of a configuration
$G^f$ is a \emph{proponent\/} of the proposal if
\emph{a majority\/} of its neighbors (including itself) are happy, namely, $f(N_v\x{N(v)})\ge 1$ (or equivalently, $\card{N_v\cap \bH}>\card{N_v\cap \bS}$); otherwise,
$v$ is an \emph{opponent\/}.
Let $P\subseteq V$ be the set of proponents of~$G^f$, and
$p=\card{P}.$ We say that $G^f$ is an \emph{approving configuration} if
$p>n/2$, namely, a majority of its vertices are proponents; otherwise, it is a \emph{disapproving configuration}.

Let $\calGnd$ denote the class of $d$-regular graphs of order $n$ with loops\footnote{Among the papers studying independent sets in regular graphs with loops we mention~\cite{GodsilNewman08}.}, i.e., graphs $G$ on $n$ vertices, each having a loop, such that $d_v=d\/$ for every $v\in V(G).$
Let $\calGndh$ denote the class of configurations $G^f$ such that $G\in\calGnd$
and it has $h$ happy vertices.
The class $\calGndh$ is called \emph{uniformly  approving}\ (respectively, \emph{uniformly  disapproving}) if it contains only approving (resp., only disapproving) configurations. $\calGndh$ is \emph{mixed} if it contains configurations of both types.
A~key question studied in this paper is the following.

\begin{problem}
\label{pm:1}
Given odd $n,d\in\N$ such that $d\le n,$ find the values of $h$ such that $\calGndh$ is uniformly approving or uniformly disapproving.
\end{problem}

Denote by $\hmind$ the minimum $h$ that ensures the existence of approving configurations for given $n$ and~$d$.
Formally, for $h=\hmind$ there exists some approving configuration in the class $\calGndh$, but for every $h<\hmind$, the class $\calGndh$ is uniformly  disapproving.
Similarly, $\hmand$ is the maximum number $h$ of happy vertices for which there exists a disapproving configuration in the class $\calGndh,$ whence for every $h\in\{\hmand+1\cdc n\},$ the class $\calGndh$ is uniformly approving.
Now Problem~\ref{pm:1} can be rewritten as follows.

\begin{manual}[problem]
{\ref{pm:1}\thmprime}\label{pm:1'}
Given odd $n,d\in\N$ such that $d\le n,$ find\x{ the value of} $\hmind$ and  $\hmand$.
\end{manual}

We also identify numbers $r$ that can be thought of as a \emph{$2$-way security gap}, in the sense that whenever the number of happy vertices deviates upwards (respectively, downwards) from $n/2$ by at least $r$, the resulting configuration class is guaranteed to be uniformly approving (resp., disapproving), and characterize the classes $\calGnd$ of regular graphs whose $2$-way security gap is sufficiently small.

Let us mention two technical features that distinguish this work from some of the previous ones.
First, we consider graphs with loops (sometimes called \emph{pseudographs}), as this provides a uniform account of the situation where the voter pays attention to their own opinion and not just to the opinions of the other neighboring voters. Such a loop is counted once in the degree of the corresponding vertex. An alternative commonly used in some of the literature is to consider local voting on the ``closed neighborhood'' $N_u\cup\{u\}$ of each vertex $u$. This is essentially equivalent,
but seems somewhat artificial and does not allow mixing voters who take and do not take into account their own opinion. Moreover, in this case, the degree of a vertex is no longer equal to the number of opinions taken into account. Considering loops is devoid of these shortcomings.

Second, we estimate the number of happy vertices that can be sufficient in some circumstances (or is definitely sufficient) for a proposal to be approved, because this keeps\x{ retains} reference to the one-quarter support level mentioned above.
Such a direct reference is lost if one estimates the difference between the number of happy and sad vertices, as is common in the literature on domination in graphs. However, for comparability, we translate our main results, Theorem~\ref{t:Classes} and Proposition~\ref{p:Foot}, into the graph domination framework in Corollaries~\ref{co:Classes} and~\ref{co:Foot} (Section~\ref{s:MajorityRegular}).

%-------------------------------------------------------------------------------------
\subsection{Basic support \lowercase{Inequalities}} %Preamble:
\label{s:Ineq}

We assume that $n$ is odd, which eliminates indefinite situations of no majority, when the numbers of proponents and opponents are equal.
Recall the \emph{handshaking lemma\/}, going back to Euler (1736).
For graphs \emph{without loops\/}, this lemma  says: \emph{The sum of vertex degrees of a graph is twice the number of the edges} (see, e.g.,~\cite{HsuLin08graph}).

In our case, since $n$ is odd, the handshaking lemma implies that $d$ (taking the loop into account) is odd too, which implies that the numbers of happy and sad neighbors of any vertex cannot be equal.
More specifically, throughout we let
\eq{e:Odd}{n ~=~ 2\kk-1,\quad d=2b-1,\quad \kk,b\in\N~.\nonumber}%-
Then a configuration $G^f$ in $\calGndh$ is an approving configuration whenever its number of proponents is at least~$\kk$, and a vertex $v$ is a proponent if and only if $N_v$ contains at least $b$ happy vertices.
For notational convenience, we view ``edge endpoints'' as distinct entities, formally defined as pairs $(v,e)$, where $v$ is a vertex and $e$ is an edge incident to~$v$. An endpoint $(v,e)$ is happy (respectively, sad) whenever $v$ is.
Therefore, in an approving configuration $G^f \in \calGndh$, the number of happy endpoints
is at least $\kk b$\x{ under the convention that a loop has one endpoint}.
On the other hand, each happy vertex contributes $d$ happy endpoints of $G^f$. Consequently, $h$ happy vertices contribute $hd$ happy endpoints in total. This number is sufficient for $G^f$ to be an approving configuration only if $hd\ge \kk b$, which implies that $\displaystyle hd \ge \frac{n+1}2\cdot\frac{d+1}2$. It follows that
\eq{e:GlobalIneq}{
\hmind ~\ge~ \frac{(n+1)(d+1)}{4d}~
}
and, expressing this bound in terms of the proportion $h/n,$
\eq{e:GlobalIneq1}{
\frac {\hmind}n ~\ge~ \frac14\big(1+n^{-1}\big)\big(1+d^{-1}\big)~.
}%+
We call \eqref{e:GlobalIneq} the \emph{global support inequality}.
\x{This inequality is necessary, but not sufficient for a configuration in $\calGndh$ to be approving.}

To obtain another condition, observe that for a configuration $G^f$ in $\calGndh$ with $h=\card{\bH}\le d/2,$ $\card{N_v\cap \bH}\le h\le d/2$ holds for every $v\in V(G^f),$ so there are no proponents in $V(G^f),$
implying that $G^f$ is a disapproving configuration.
{Therefore,
$\displaystyle h ~\ge~ \frac{d+1}2$,
is also necessary for a configuration in $\calGndh$ to be approving, implying that}
\eq{e:LocalIneq}{\hmind ~\ge~ \frac{d+1}2~.}
We refer to this inequality as the \emph{local support inequality}.

{Combining the global and local support inequalities} \eqref{e:GlobalIneq} and \eqref{e:LocalIneq}, we have the following.
\begin{proposition}
\label{p:Ineq}
$\displaystyle \hmind\,\ge\,{\frac12(d+1)}\x{\frac{d+1}2} \max\left\{\frac{n+1}{2d};\, 1\right\}.$
\end{proposition}

\begin{figure}[ht] %\x{ht}
\centerline{\small (a)\hspace{25em}(b)}
\centerline{
\includegraphics[height=5.5cm]{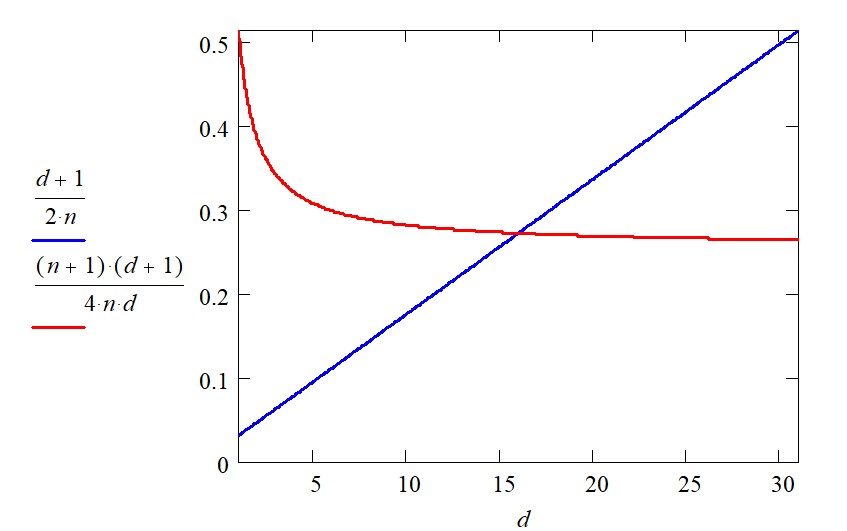}\;\; %[width=8cm]
\includegraphics[height=5.5cm]{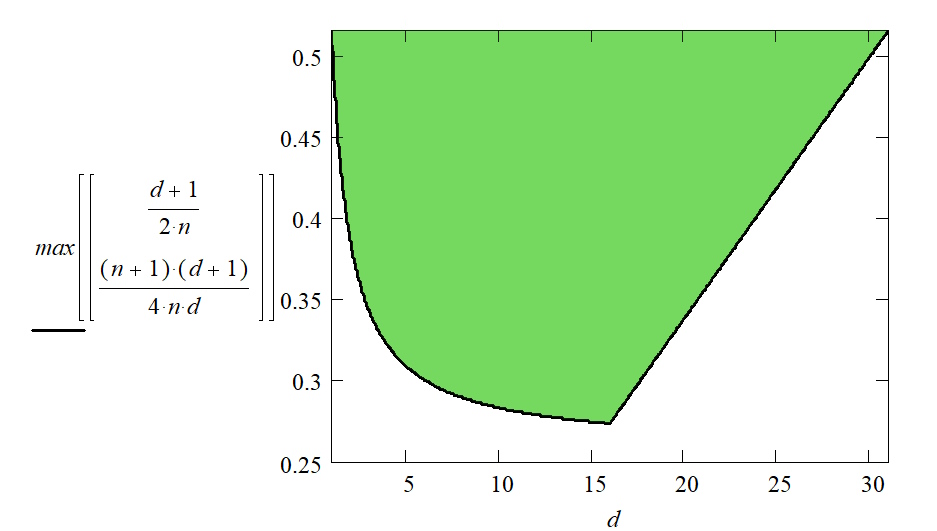}
}
\caption{Boundary curves of the global and local support inequalities~(a) and the support domain~(b) for~$n=31$ in the coordinates $d$ and~$h/n$ \label{f:Curve}}
\end{figure}

Conditions \eqref{e:GlobalIneq} and \eqref{e:LocalIneq} complement each other. Let us illustrate this by considering networks of order $n=31.$ For such networks, the corresponding boundary curves and the intersection of the support domains in the coordinates $d$ and $h/n$ are shown in Fig.~\ref{f:Curve}.

%-------------------------------------------------------------------------------------
%\setcounter{figure}{-1}
%-------------------------------------------------------------------------------------

The intersection point of the boundary curves of the local and global support inequalities is
\eq{e:Intersect}{\check d(n)~=~\frac{n+1}2~=~\kk;\quad\quad \check h(n)~=~\frac{\check d(n)+1}2~=~\frac{n+3}4~.} %=b%+
Indeed, to find this point, it suffices to solve the equation
$$
\frac{n+1}2\cdot\frac{d+1}{2d} ~=~ \frac{d+1}2~.
$$
Since $\frac14(n+1)(1+d^{-1})$ (respectively, $\frac{d+1}2$) is decreasing (resp., increasing) in $d,$
any approving configuration satisfies
\eq{e:minh}{h ~\ge~ \check h(n) ~=~ \frac{n+3}4~.}%+
Note that if $\check d(n)$ is even, then no graph matches the intersection point.

In what follows, we call graphs in $\calGnd$ with $d\ge\frac{n+1}2$ \emph{high-degree $d$-regular graphs} and graphs with $d<\frac{n+1}2$ \emph{low-degree $d$-regular graphs}.
Observe that for high-degree graphs,
the local support inequality is no weaker than the global support inequality, while in contrast, for low-degree graphs,
the global support inequality is stronger.
In Section~\ref{s:Constru}, we construct high-degree and low-degree regular approving configurations with any fixed $n$ and $d$ and find $\hmind$ and $\hmand$ for these $n$ and~$d.$

%-------------------------------------------------------------------------------------
\section{Characterization of configuration classes}
\label{s:Constru}

{Is Proposition} \ref{p:Ineq},
the conjunction of support inequalities \eqref{e:GlobalIneq} and \eqref{e:LocalIneq}, tight? That is, does it provide a sufficient condition for a class $\calGndh$ to contain an approving configuration? We show that this is true with two exceptions.
{The following two subsections carry out this analysis separately for \emph{high-degree} configurations (defined as ones where $d\ge (n+1)/2$) and \emph{low-degree} configurations (where $d<(n+1)/2$).}

%-------------------------------------------------------------------------------------
\subsection{High-degree approving configurations with minimum number of happy vertices}
\label{ss:HighConf}

We first consider high-degree graph classes $\calGnd,$ that is, where
\eq{e:HighDegree}{d ~\ge~ \frac{n+1}2~.}

Our goal in this section is to investigate whether for all high-degree configurations,
the local support inequality \eqref{e:LocalIneq} is tight, i.e., $\hmind=(d+1)/2,$ or more explicitly,
there are approving configurations in
$\calGndh$ with $d=2h-1$.

{In the extremal case of $n=d$,} $\calGnd$ contains only the complete graph, and any corresponding configuration with $h=(d+1)/2$ is an approving configuration, since every vertex of it is a proponent.
Thus, the lower bound determined by the local support inequality \eqref{e:LocalIneq} for {odd} $d=n$ is {tight, i.e., $\hmind=(d+1)/2$}.

Next, consider configurations in $\calGndh$
with $n>d$.
For any such configuration with odd $n,$ there exists some nonnegative integer $t$ such that $n=2h+1+2t,$
so the configuration belongs to class $\calG_{2h+1+2t|2h-1|h}$ with
\eq{e:HighDegree_t}{0 ~\le~ t ~\le~ h-2\,,}%+
for which
the inequality $t\le h-2$ is equivalent to the condition \eqref{e:HighDegree} of high degree.

The integer range $[0,h-2]$ defined for $t$ in Eq.  \eqref{e:HighDegree_t}
can be split into three parts, some of which can be empty,
\eq{e:3cases}{
[0;\;h-2] ~=~ \left[0;\,\left\lfloor\tfrac{h-3}2\right\rfloor\right]
~\cup~ \left[\left\lfloor\tfrac{h-3}2\right\rfloor+1;\;h-3\right]
~\cup~ \{h-2\},
}%+
with different regularities, which will be considered separately.
\medskip

We use the following notation.

\begin{itemize}
\item
$H_w$ is any configuration of order $w$ with $w$ happy vertices;
\item
$S_w$ is any configuration of order $w$ with $w$ sad vertices;
\item When $S_w$ is a subconfiguration of a configuration $G^f$ in which all vertices of $S_w$ are proponents, we sometimes denote $S_w$ by~$PS_w.$
\item $C_w$ is any configuration whose graph is a circulant graph of order~$w$ with loops.
\end{itemize}

For our constructions, we use a number of graph composition operations,
which are defined as follows.

\begin{itemize}
\item $G=\conne{G'}{G''}$ is the \emph{join} of the graphs $G'$ and $G''$,
such that $V(G')\cap V(G'')=\varnothing$:
$V(G)=V(G')\cup V(G'')$ and $E(G)=E(G')\cup E(G'')\cup (V(G')\times V(G'')).$
\item
The graph $G=\connee{G'}{G''}{G'''}$ is {the \emph{double join} of the graphs $G'$, $G''$ and $G'''$,}
defined similarly with $E(G)=E(G')\cup E(G'')\cup E(G''')\cup (V(G')\times V(G''))\cup (V(G'')\times V(G''')).$
\item $G=\doubling{G'}{G''}$ is {the \emph{$2$-to-$1$ matching} of the graphs $G'$ and $G''$,} where $V(G)=V(G')\cup V(G'')$ is a partition of $V(G),$ $\card{V(G')}=2\card{V(G'')},$ every vertex of $G'$ has exactly one neighbor in $V(G''),$ and every vertex of $G''$ has exactly two neighbors in $V(G').$
\item $G=\mediastinum{G'}{G''}{G'''}$ is {the \emph{hybrid join} of the graphs $G'$, $G''$ and $G'''$,} where $G'$ is connected to $G''$ by the ${\succ\hspace{-5pt}-}$ %>
operation, while $G''$ is connected to $G'''$ by the ${\equiv\!\equiv}$ operation.
\item $G=[G'\,\;G'']$ is the \emph{union} of graphs $G'$ and $G''$, such that $V(G')\cap V(G'')=\varnothing$:\x{a graph, where} $V(G)=V(G')\cup V(G'')$\x{ is a partition of $V(G)$} and $E(G)=E(G')\cup E(G'').$
\end{itemize}

All of these operations apply also to binary labeled graphs, i.e., configurations, where the function $f$ carries over to $V(G)$, so $\bH(G)=\bH(G')\cup \bH(G'')$ or $\bH(G)=\bH(G')\cup \bH(G'') \cup \bH(G''')$, as the case may be.
\medskip

We make frequent use of the following well-known lemma.

\begin{lemma}
\label{lem:reg-graph}
For any $n,d\in\N,$ a $d$-regular $n$-vertex graph with loops exists if and only if $d\le n$ and\/ $(d+1)n$ is even.
\end{lemma}

In Lemma~\ref{lem:reg-graph}, necessity follows from the handshaking lemma, and sufficiency follows from the possibility of constructing a circulant graph with appropriate parameters $n$ and~$d$ (see, e.g., \cite[p.\,11]{Melnikov98exercises}).

\medskip
Having \eqref{e:3cases} in mind, let us consider three cases.

\medskip
{\bf Case 1.} Suppose that
\eq{e:Case1}{0 ~\le~ t ~\le~ \frac{h-3}2~.}%+

Since $n=2h+1+2t$ and $d=2h-1,$ segment \eqref{e:Case1} corresponds to $d\in\left\{2\!\left\lceil\frac {n+2}3\right\rceil-1\cdc n-2\right\}$. It degenerates when $h\le2$ or $n\le5$.
%\aB{Since $n=2h+1+2t$ and $d=2h-1,$ segment \eqref{e:Case1} corresponds to $d\in\left\{2\!\left\lceil\frac {n+2}3\right\rceil-1\cdc n-2\right\}$. It degenerates when $h\le2$ or $n\le5$}.

\begin{proposition}\label{p:High1}
Under Ineq.
\eqref{e:Case1}$,$ there exists an approving configuration $G^f(1)\in\calG_{2h+1+2t|2h-1|h}$ of the form $G^f(1)=\composiTwo{h+2t+1}{h}$.
\end{proposition}

\begin{proof}
Define $G^1=(V^1,E^1)$ and $G^2=(V^2,E^2)$ to be regular graphs with loops such that $\card{V^1}=h$, $\card{V^2}=h+2t+1$ and $V^1\cap V^2=\varnothing$.
The degrees of $G^1$ and $G^2$ are $(2h-1)-(h+2t+1)=h-2t-2$ and $(2h-1)-h=h-1,$ respectively. Note that $h-2t-2\ge1$ since $t\le(h-3)/2$; $h-1>1$ since $h>2$. A possible choice of $G^1$ and $G^2$ are appropriate circulant graphs with loops. These exist due to Lemma~\ref{lem:reg-graph}, because $(h-2t-1)h$ and $h(h+2t+1)$ are even.
Set $f^1(v)=1$ for all $v\in V^1$ and $f^2(v)=-1$ for all $v\in V^2.$ Now define $H_h=(V^1,E^1,f^1)$ and $S_{h+2t+1}=(V^2,E^2,f^2)$.

With this definition and $f=f^1\cup f^2$ the number of happy vertices in $G^f(1)=\composiTwo{h+2t+1}{h}$ is $h$, its order is $(h+2t+1)+h=2h+1+2t$, and every vertex in it has degree $2h-1$. Hence $G^f(1)\in\calG_{2h+1+2t|2h-1|h}.$

$G^f(1)$ is an approving configuration since
 (i)~every vertex of $S_{h+2t+1}$ is a proponent (as $h$ out of $2h-1$ neighbors of it are happy), so that $S_{h+2t+1}=PS_{h+2t+1}$ and (ii)~$\card{V^2}>\card{V^1}.$
\end{proof}

Note that in the configuration $G^f(1)$ used in the proof of Proposition~\ref{p:High1}, the set of proponents coincides with the set of sad vertices.
\medskip

{\bf Case 2.} Suppose now that
\eq{e:Case2}{\frac{h-3}2 ~<~ t ~\le~ h-3\,.}%+

Then $d\in\left\{2\!\left\lceil\frac{n+1}4\right\rceil+1\cdc2\!\left\lceil\frac {n+2}3\right\rceil-3\right\}$. This range is empty when $h\le3$ or $n\in\{1,3,5,7,9,13\}$.
%\aB{Then $d\in\left\{2\!\left\lceil\frac{n+1}4\right\rceil+1\cdc2\!\left\lceil\frac {n+2}3\right\rceil-3\right\}$. This range is empty when $h\le3$ or $n\in\{1,3,5,7,9,13\}$.}
%
%$n\le5$ or $n=9$.}

\begin{proposition}\label{p:High2}
Under
Ineq.\,\eqref{e:Case2}$,$ there exists an approving configuration $G^f(2)\in\calG_{2h+1+2t|2h-1|h}$ of the form $G^f(2)=\composiThree{2t-h+3}{2h-2}{h}$.
\end{proposition}

\begin{proof}
Start by defining three graphs, $G^1=(V^1,E^1)$, $G^2=(V^2,E^2),$ and $G^3=(V^3,E^3)$, as follows.
Select disjoint vertex sets $V^1$, $V^2,$ and $V^3$ with $\card{V^1}=h$, $\card{V^2}=2h-2,$ and $\card{V^3}=2t-h+3.$
Note that $2h-2>4$ as Case~2 occurs only when $h>3$ and $2t-h+3>0$ as $t>(h-3)/2$ by \eqref{e:Case2}. Let the edge sets $E^1$ and $E^3$ contain only loops, and let $G^2$
be a regular graph of degree $(2h-1)-(2t-h+3)-h=2(h-t-2)$ with loops. We have $2\le2(h-t-2)\le h-1$ by~\eqref{e:Case2}.
A possible choice for $G^2$ is an appropriate circulant graph; having an even number of vertices, $G^2$ does not violate the conditions of Lemma~\ref{lem:reg-graph}.
Set $f^1(v)=1$ for all $v\in V^1$, and set $f^2(v)\equiv-1$ and  $f^3(v)\equiv-1$ on their entire domains.
Now set $H_h=(V^1,E^1,f^1),$ $S_{2h-2}=(V^2,E^2,f^2),$ and $S_{2t-h+3}=(V^3,E^3,f^3).$

With this definition and $f=f^1\cup f^2\cup f^3,$ the number of happy vertices in $G^f(2)=\composiThree{2t-h+3}{2h-2}{h}$ is $h$, its order is $(2t-h+3)+(2h-2)+h=2h+1+2t$, and every vertex in it has degree $2h-1$ by construction. Hence $G^f(2)\in\calG_{2h+1+2t|2h-1|h}.$

$G^f(2)$ is an approving configuration, since
(i)~every vertex of $S_{2h-2}$ is a proponent ($h$ out of $2h-1$ neighbors of it are happy), so that $S_{2h-2}=PS_{2h-2}$ and $p=2h-2;$
 (ii)~by \eqref{e:Case2}, $n=2h+1+2t\le2h+1+(2h-6)<2(2h-2)=2p.$
\end{proof}

Observe that in the configuration $G^f(2)$ constructed in the proof of Proposition~\ref{p:High2}, the set of proponents is a proper subset of the set of sad vertices. The noted properties of $G^f(1)$ and $G^f(2)$ contrast with the setting of Fishburn et al.~\cite{Fishburn86local}, presented in Subsection~\ref{ss:Related},
where only happy vertices may become proponents, and moreover, all of them must be proponents for the approval of a proposal.\!
\medskip

{\bf Case 3.} $t=h-2.$
\medskip

In this case
%(which does not apply to $h=1$ by virtue of \eqref{e:HighDegree_t}),
%
(which degenerates when $h=1$ or $\frac{n+1}4\in\N$), $d=\frac{n+1}2$ and
%
%\aB{(which degenerates when $h=1$ or $\frac{n+1}4\in\N$), $d=\frac{n+1}2$ and}
%
$\calG_{2h+1+2t|2h-1|h}=\calG_{4h-3|2h-1|h}.$

\begin{proposition}\label{p:High3}
The class $\calG_{4h-3|2h-1|h}$ is uniformly  disapproving $($i.e.$,$ it contains no approving configuration$)$ when $h\in\{2,3\}.$ If $h>3,$ then there exists an approving configuration $G^f(3)\in\calG_{4h-3|2h-1|h}$ of the form $G^f(3)=\mediastinum{S_{2h-2}}{S_{h-1}}{H_h}$.
\end{proposition}

\begin{proof}
Since $d=\frac{n+1}2$ and $h=\frac{d+1}2,$ Case~3 refers to the intersection point \eqref{e:Intersect} of the boundary curves of the local and global support inequalities. As the global support \emph{equality\/} holds at this point, only proponents are adjacent to happy vertices. Therefore, because of the loops, all happy vertices are proponents.

Consider the case $h=2.$ If $G^f \in \calG_{4h-3|2h-1|h}=\calG_{5|3|2}$ is an approving configuration, then it has three proponents, which are the two happy vertices and one sad vertex. They must be adjacent to each other, and as $d=3,$ they are not adjacent to the two remaining vertices. This leaves no chance for the two remaining vertices to have degree~3. This contradiction indicates that no configuration in $\calG_{5|3|2}$ is an approving configuration.

Now consider $h=3.$ If $G^f \in \calG_{4h-3|2h-1|h}=\calG_{9|5|3}$ is an approving configuration, then there are five proponents, which are the three happy and two sad vertices. The four remaining vertices cannot all have degree 5 because each of them has at most four connections apart from the five proponents, whence they must have at least four connections to the five proponents altogether, while the aforementioned proponents can provide only two. This contradiction shows that no configuration in $\calG_{9|5|3}$ can be an approving configuration.

Finally suppose $h>3.$ We prove that there exists an approving configuration $G^f(3)\in\calG_{4h-3|2h-1|h}$ of the form $G^f(3)=\mediastinum{S_{2h-2}}{S_{h-1}}{H_h}.$
Define the adjacency within the configurations $S_{2h-2},$ $S_{h-1},$ and $H_h$ as follows. $S_{2h-2}$ and $H_h$ are based on complete graphs with loops; each vertex of $S_{h-1}$ is adjacent to all vertices of $S_{h-1}$ (including itself), except for two, which satisfies the conditions of Lemma~\ref{lem:reg-graph}
for all $h>3,$ since $(h-1)(h-2)$ is even.
To obtain such a graph, we can construct an appropriate circulant graph. The opinion function of $G^f(3)$ assigns $1$ to the vertices of $H_h$ and $-1$ to the remaining vertices.\x{ of $S_{2h-2}$ and  $S_{h-1}.$}

%3.2.
With this definition, the order of $G^f(3)=\mediastinum{S_{2h-2}}{S_{h-1}}{H_h}$ is $4h-3$, the number of happy vertices is $h$, and the degree of each vertex is $2h-1$ (specifically, the degree of the vertices of $H_h,$ $S_{2h-2},$ and $S_{h-1}$ are $(h-1)+h,$ $(2h-2)+1,$ and $2+(h-1-2)+h,$ respectively). Therefore, $\mediastinum{S_{2h-2}}{S_{h-1}}{H_h}\in\calG_{4h-3|2h-1|h}.$

%3.3.
The defined configuration $G^f(3)=\mediastinum{S_{2h-2}}{S_{h-1}}{H_h}$ is an approving configuration, since all vertices of $H_h$ and $S_{h-1}$ are proponents as they are adjacent to $h=\frac{d+1}2$ happy vertices, while the number of these proponents is $2h-1=\frac{n+1}2.$
\end{proof}

Now let us\x{ specify} find the minimum values of $h$ for which the classes $\calG_{5|3|h}$ or $\calG_{9|5|h}$ contain at least one approving configuration.

\begin{proposition}\label{p:High4}
There exist approving configurations in the classes $\calG_{5|3|3}$ and\/ $\calG_{9|5|4}.$
\end{proposition}

\begin{proof}
An approving configuration in the class $\calG_{5|3|3}$ is the 3-regular cyclic graph with loops on 5 vertices with an arbitrary location of 3 happy vertices.

An approving configuration in $\calG_{9|5|4}$ is shown in Fig.~\ref{f:G954}.
%----------------------------------------------------
\begin{figure}[ht] %\x{ht}
\centerline{\small \!\!\!\!\!\!(a)\hspace{19.3em}(b)}
\medskip
\centerline{
\includegraphics[width=12em]{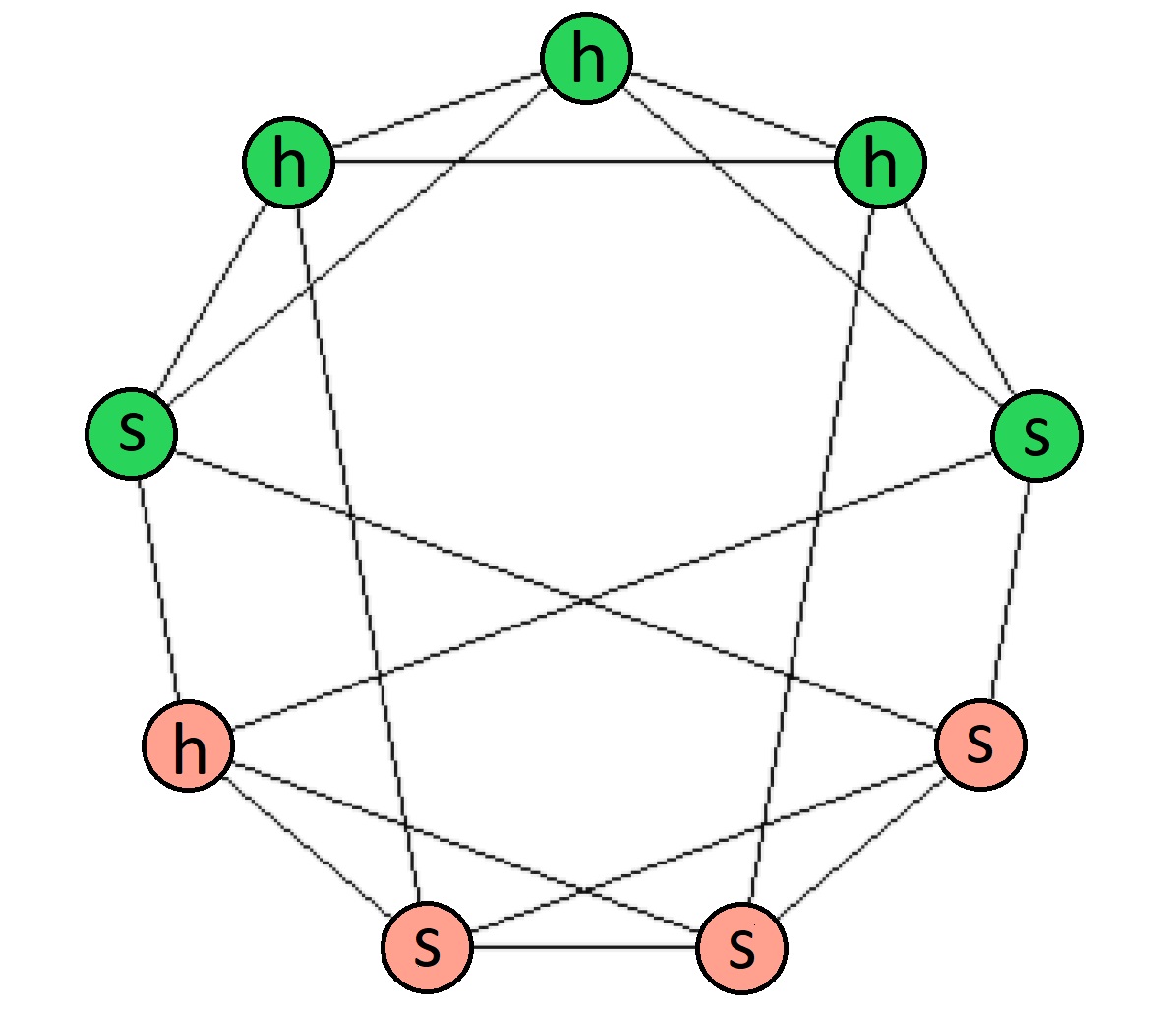}
\qquad \qquad \qquad
\includegraphics[width=14em]{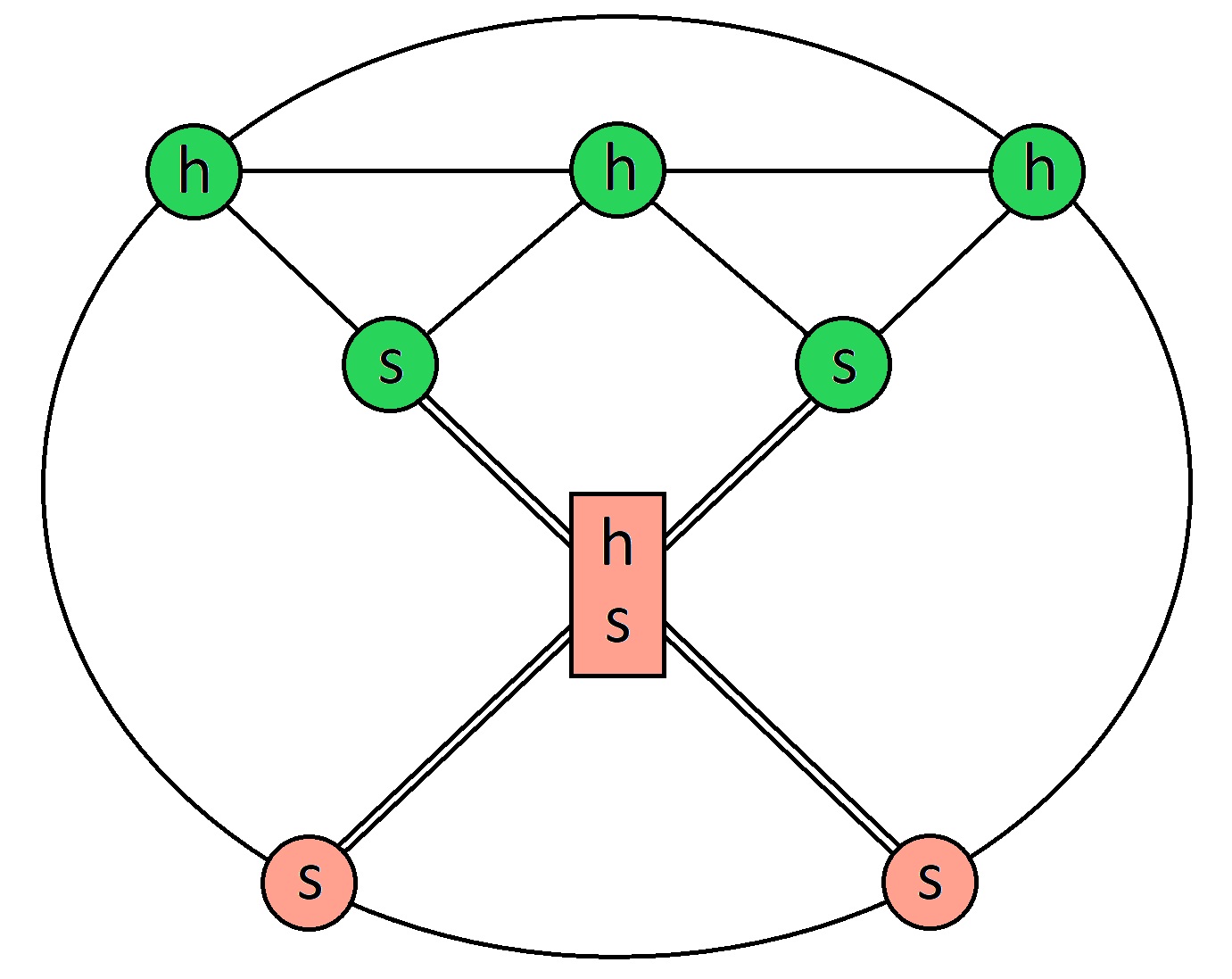}
}
\caption{An approving configuration $G_{9|5|4}$: its standard (a) and quasi-plane (b) drawings. Happy and sad vertices are marked by h and s, respectively. In (b), the vertices $h$ and $s$ in the rectangle are both adjacent to the four vertices connected to the rectangle by double edges. Loops are not shown. The green (respectively, red) vertices are the proponents (resp., opponents).}\label{f:G954}
\end{figure}
%----------------------------------------------------
This configuration is obtained from one of sixteen 4-regular connected graphs of order 9 without loops by adding loops for each vertex; the loops are not shown in the figure. The graph has five proponents, which are shown in green, hence it is an approving configuration.
\end{proof}

Proposition~\ref{p:High4} completes the construction of high-degree approving configurations with a minimum number of happy vertices. Combined with Proposition~\ref{p:Ineq}, we get the following lemma, which summarizes the bounds on $\hmind$ for high-degree configurations.

\begin{lemma}
\label{lem:high-deg characterization}
For high degree configuration classes$,$ $\hmind=(d+1)/2,$ with two exceptions\/$:$ $\hmin(5,3)=3$\/ and\/ $\hmin(9,5)=4$.
\end{lemma}

{It is convenient to represent this (and subsequent) results using indicator function notation. For a universal set $U$, a subset $A\subseteq U$
%$A\subseteq  \{(n,d) \mid n,d \mbox{ odd},~ d\le n\}$,
and an element $x\in U$, let ${\bm1}_A(x)$ be\/ $1$ if $x\in A$ and\/ $0$ otherwise.}

%Pavel{In Section 5 we also use $\bm1_A(x)$ with a scalar argument $x=n.$ So it would be good to give ${\bm1}_A(x)$ a universal definition.}
%David: ok.

%\vspace{-.35em}
\begin{corollary}\label{c:high-deg hmin}
For high degree configuration classes\/$,$
$$
\hmind=\x{\frac12(d+1)}\frac{d+1}2+\indi{(5,3),(9,5)}{(n,d)}\,.
$$
%where ${\bm1}_A(x)$ is the indicator function of a subset $A$ of the universal set to which $x$ belongs\/$,$ namely\/$,$ $\bm1_A(x)$ is\/ $1$ if $x\in A$ and\/ $0$ otherwise.
\end{corollary}
%David{This definition might still be unclear. Perhaps we should place the definition before the corollary, and define the relevant universal set here, as $U=\{(n,d) \mid n,d \mbox{ odd }, (n+1)/2 \le d\le n\}$.}
%Pavel{Another option is to add it to the bulleted list on page 10.}
%David{That list defines the graph components used in our constructions, so it would be odd to place it there. I would prefer defining it close to where we use it for the first time.}

%----------------------------------------------------
\subsection{Low-degree approving configurations with minimum number of happy vertices}
\label{ss:LowConf}

Let us\x{ now} turn to low-degree graph classes $\calGnd,$ i.e.,
%where
classes with $d<(n+1)/2=\kk$
%\eq{e:LowDegree}{d~<~\frac{n+1}2~=~\kk.}
(in other form,  $n>1$, $d\in\left\{1\cdc2\!\left\lfloor\frac{n+1}4\right\rfloor-1\right\}$).
%
%\aB{(in other form,  $n>1$, $d\in\left\{1\cdc2\!\left\lfloor\frac{n+1}4\right\rfloor-1\right\}$)}.
We show that for such configurations,
the global support inequality \eqref{e:GlobalIneq} is tight, or more explicitly, for all $n,$ $d,$ and $h$ satisfying \eqref{e:GlobalIneq}, there exist approving configurations in $\calGndh$.

Since $n>2d-1$ and $n$ is odd, there exists $t\in\{0,1,2,\ldots\}$ such that $n=2d+1+2t.$
Therefore, we are dealing with configurations from the class $\calG_{2d+1+2t|d|h}.$

\begin{proposition}\label{p:Low}
For any class $\calG_{2d+1+2t|d|h}$ satisfying the global support inequality~\eqref{e:GlobalIneq}$,$
there exists an approving configuration $G^f(4)\in\calG_{2d+1+2t|d|h}$ of the form $G^f(4)=[S_{d+t}\;\,C_{d+t+1}]$.
%Each class $\calG_{2d+1+2t|d|h}$ satisfying the global support inequality~\eqref{e:GlobalIneq} contains an approving configuration.
\end{proposition}

\begin{proof}
For any class $\calG_{2d+1+2t|d|h}$ consider disconnected configurations of the form $G^f(4)=[S_{d+t}\;\,C_{d+t+1}]$ with components $S_{d+t}$ and $C_{d+t+1}$ of order $d+t$ and $d+t+1,$ respectively, where the circulant component $C_{d+t+1}$ will be defined later, while $S_{d+t}$ is $d$-regular with loops and all of its vertices are sad. The latter is possible due to Lemma~\ref{lem:reg-graph}, since $d$ is odd.
Before\x{ Prior to} defining the structure of $C_{d+t+1}=C_{\kk},$\x{ of order $d+t+1,$} we estimate the ratio of the number of happy and sad vertices in it provided that %$=[S_{d+t}\;\,C_{d+t+1}]$
$G^f(4)\in\calG_{2d+1+2t|d|h}$ and \eqref{e:GlobalIneq} is satisfied.

Substituting $n=2\kk-1$ in \eqref{e:GlobalIneq} yields
\eq{e:k_Low}{\kk ~\le~ \frac{2h}{1+d^{-1}}~.}%+

Denote by $s'$ the number of sad vertices in $C_{d+t+1}.$ It follows from \eqref{e:k_Low} that
\eq{e:s'}{s'~=~\kk-h ~\le~ \frac{d-1}{d+1}h ~<~ h~.}%+

The construction of $C_{d+t+1}$ starts with indexing all the sad vertices by numbers $2, 4\cdc 2s'$ and $s'$ happy vertices by $1, 3\cdc 2s'-1.$

We now estimate the average distance $\rho$ between the indices of the remaining subsequent happy vertices if we insert them into the above sequence freeing up one index for each and considering this sequence ${\!\!\!}\pmod \kk$, i.e., as a cyclic one. Since the number of the remaining happy vertices is $h-s'$ and $\kk=h+s',$ using \eqref{e:s'} we obtain
\eq{e:rho}{\rho~=~\frac{h+s'}{h-s'}~\le~ d~.\nonumber}%-

Consequently, the remaining happy vertices can be indexed in such a way that the cyclic ${\!\!\!}\pmod \kk$ distance between the subsequent added happy vertices does not exceed~$d.$ As a result, any cyclic segment of length $d-1$ of this sequence contains at least one added happy vertex and therefore (recalling that $d$ is odd) it contains more happy vertices than sad vertices.

We now complete the construction of $C_{\kk}=C_{d+t+1}$ by making each vertex adjacent to itself, to $\frac{d-1}2$ consecutive vertices to the left, and $\frac{d-1}2$ consecutive vertices to the right of it in the above circle of indices (obtaining a Harary circulant graph with loops). This can be done, since for low-degree configurations
%\eqref{e:LowDegree}
$d<\kk.$

By construction, %$[S_{d+t}\;\,C_{d+t+1}]$
$G^f(4)\in\calG_{2d+1+2t|d|h}.$ Since for every vertex $v$ of $C_{d+t+1}$ its set of neighbors %neighborhood
$N_v$ forms a length $d-1$ cycle's segment containing more happy vertices than sad vertices, every $v$ is a proponent. The number of such proponents is $\kk=\frac{n+1}2$ and thus, the constructed configuration $[S_{d+t}\;\,C_{d+t+1}]$ is approving.
\end{proof}

Combined with Proposition~\ref{p:Ineq}, we get the following lemma, which summarizes the bounds on $\hmind$ for low-degree configurations.

\begin{lemma}
\label{lem:low-deg characterization}
For low degree configuration classes$,$
$$\hmind ~=~ \left\lceil\frac{(n+1)(d+1)}{4d}\right\rceil.$$
\end{lemma}

%-------------------------------------------------------------------------------------
%\section{Characterization of configuration classes}
%\label{s:WhichCan}

%-------------------------------------------------------------------------------------
\subsection{The complete classification}
%{Classes that contain approving configurations}
\label{s:Classes}

Combining the results of the last two subsections, we completely characterize the classes $\calGndh$ of regular configurations that contain approving configurations, as well as the uniformly approving or disapproving configuration classes.

\begin{proposition}\label{p:Canapproval}
For odd\/ $n,d\in\N$ such that $d\le n,$
$$
\hmind~=~
\left\lceil\frac{d+1}{2}\cdot\max\left\{\frac{n+1}{2d}\,,\,1\right\}\right\rceil+\indi{(5,3),(9,5)}{(n,d)}\,.
$$
\vspace{-1.2em}
\end{proposition}

\begin{proof}
This follows from Corollary~\ref{c:high-deg hmin}, Lemma~\ref{lem:low-deg characterization}, and the fact that $\max\{\frac{n+1}{2d},\,1\}$ is $\frac{n+1}{2d}$ in the case of low degree configuration classes and $1$ in the opposite case.
\end{proof}

\begin{proposition}\label{p:Conforming}
For odd\/ $n,d\in\N$ such that $d\le n,$
$\displaystyle
\hmand\,=\,n-\hmind.$
\end{proposition}

\begin{proof}
For any configuration $G^f \in \calGndh$ consider the configuration $G^{-f} \in \calG_{n|d|n-h}$
obtained from $G^f$ by replacing $f$ with $-f,$\x{ multiplying $f$ by~$-1,$} which is equivalent to interchanging the sets $\bH(G^f)$ and $\bS(G^f)=V(G^f)\setminus \bH(G^f)$.
Observe that $G^f$ is an approving configuration if and only if $G^{-f}$ is not. Consequently, a class $\calGndh$
is uniformly approving if and only if the class $\calG_{n|d|n-h}$
is uniformly disapproving, and hence
{the thresholds $\hmind$ and $\hmand$ are connected as in the proposition.}
\x{$\hmind$ that isolates uniformly disapproving classes is connected with $\hmand$ isolating uniformly approving classes by the equation $\hmand\,=\,n-\hmind$.}
\end{proof}

Combining Propositions~\ref{p:Canapproval} and~\ref{p:Conforming} we obtain the following theorem, which characterizes approving, mixed, and disapproving classes $\calGndh$ and resolves Problem~\ref{pm:1} formulated in Section~\ref{s:Problem}.

\begin{theorem}\label{t:Classes}
For odd\/ $n,d\in\N$ such that $d\le n$ and $h\in\N$,
a class\/ $\calGndh$ of binary labeled $d$-regular graphs of order $n$ with $h$ happy vertices is\/$:$

\medskip
\begin{tabular}{llrcl}
$-$ uniformly  disapproving &\! iff &\!$0\,\le\,$\lI&$h$&\lI$\,<\,\hmind\;;$
\\
$-$ mixed       &\! iff\: &\!$\hmind\,\le\,$\lI&$h$&\lI$\,\le\,\hmand\lB{\mathstrut}\,;$
\\
$-$ uniformly  approving    &\! iff &\!$\hmand\,<\,$\lI&$h$&\lI$\,\le\, n\,,$
\end{tabular}

\medskip
\noindent
where
\eq{e:hmimd}{
\hmind~=~
\left\lceil\frac{d+1}{2}\cdot\max\left\{\frac{n+1}{2d}\,,\,1\right\}\right\rceil+\indi{(5,3),(9,5)}{(n,d)}
}
and\, $\hmand=n-\hmind\,.$
\end{theorem}

\x{The bounds $\beta_1(n,d)$ and $\beta_2(n,d)$ of the region of mixed classes
on the plane with coordinates $d$ and $h$ are symmetric w.r.t. the line $h=\frac n2.$}
\begin{corollary}\label{co:trusting}
Under the assumptions of Theorem~$\ref{t:Classes},$ \
$(\hmind-1)/n$ is the greatest number of the form $\alpha_1=l/n$ $(l\in\N)$ such that the two-stage majority procedure on $d$-regular graphs of order $n$ is $\alpha_1$-protecting\/$;$
$\hmand/n$ is the smallest number $\alpha_2$ such that this procedure is $\alpha_2$-trusting.
\end{corollary}

\begin{proof}
This follows from Theorem~\ref{t:Classes} and the relevant definitions.
\end{proof}

\begin{figure}[ht] %\x{ht}
\centerline{\small (a)\hspace{25em}(b)}
\centerline{
\includegraphics[height=15.0em]{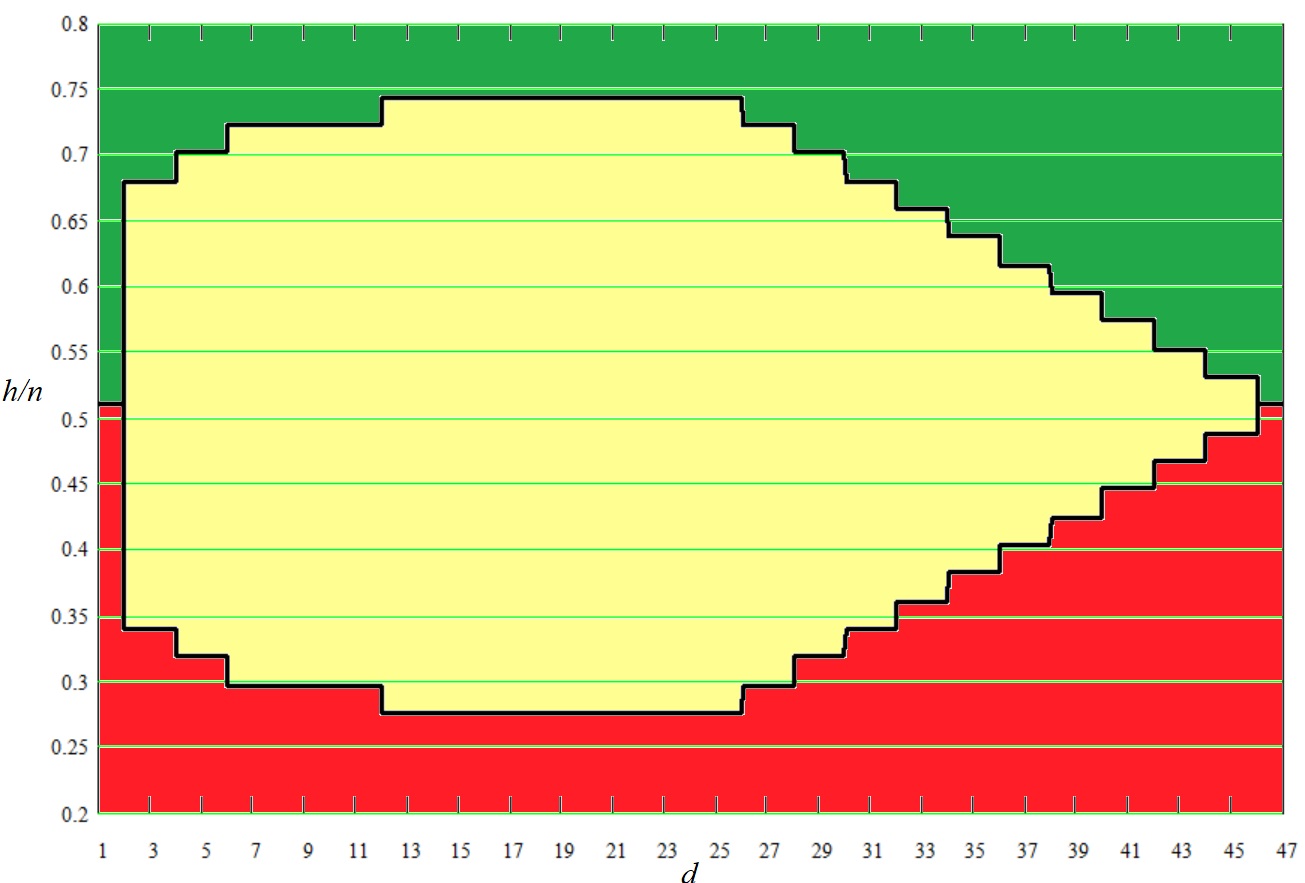}\;\; %[width=23.0em]
\includegraphics[height=15.0em]{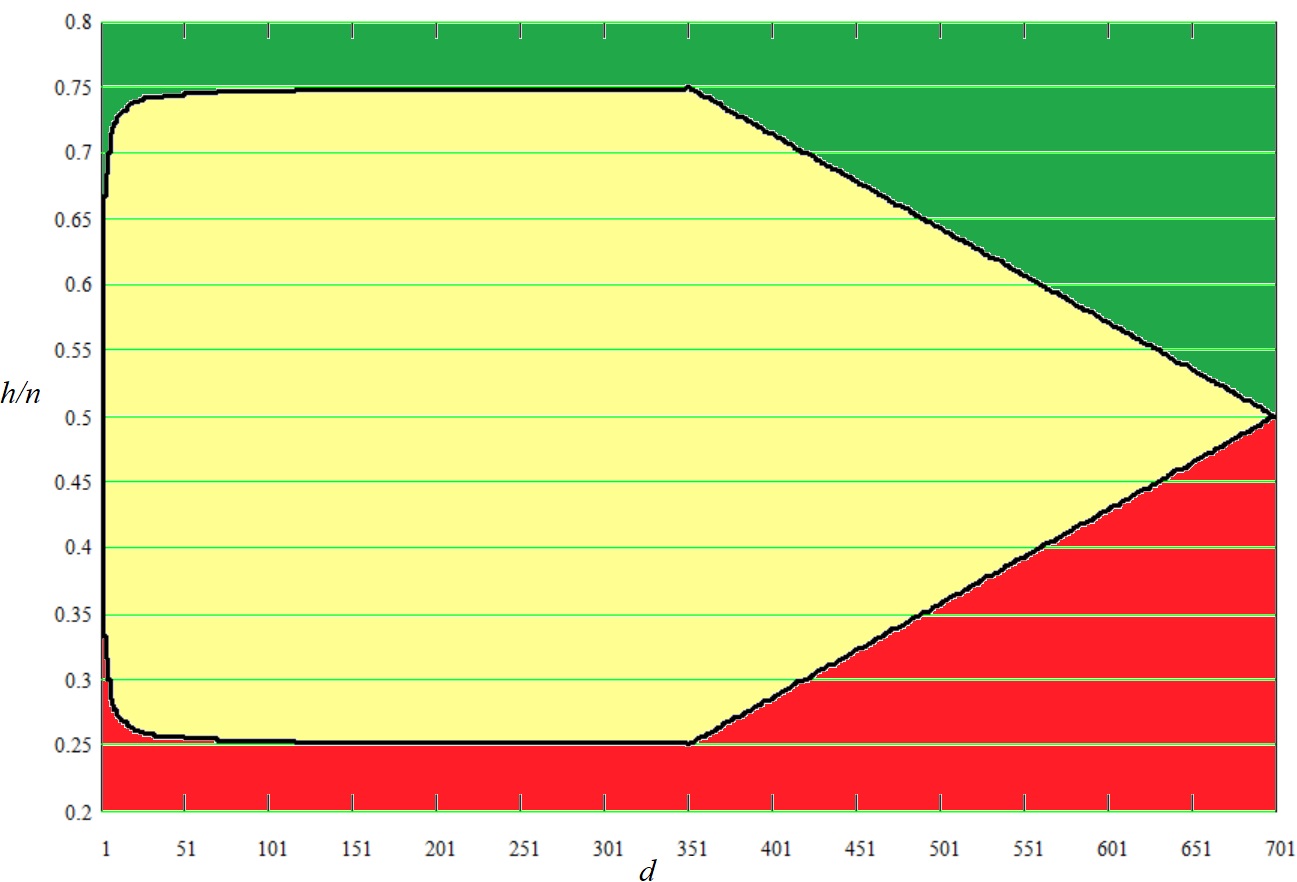} %
}
\caption{Regions of mixed (yellow), uniformly  disapproving (red), and uniformly  approving (green) classes $\calGndh$ of binary labeled $d$-regular graphs of order $n$ with $h$ happy vertices in the coordinates $d$ and $h/n$ for (a)~$n=47$ and (b)~$n=701$.
\label{f:ThreeDom}}
\end{figure}

Theorem~\ref{t:Classes} is illustrated by Fig.~\ref{f:ThreeDom}, where the regions of mixed (yellow), disapproving (red), and approving (green) classes in the coordinates $d$ and $h/n$ are shown for $n=47$ and $n=701$ taking into account that $h$ is integer and $d$ is odd. More precisely, the upper separating boundary contains the minimum values of $h/n$ that belong to the domain of \uniform approval, while the lower separating boundary contains the minimum values of $h/n$ that \emph{do not belong\/} to the domain of \uniform disapproval.

Note that these boundaries are not defined symmetrically w.r.t. $h=\frac n2,$ because any integer-valued symmetric boundaries would give either overlapping regions (since the region of mixed classes is empty at $d\in \{1,n\}$) or a gap between them. It can be easily demonstrated that the aforementioned boundaries shown in Fig.~\ref{f:ThreeDom} are symmetric w.r.t. the line $h=\frac {n+1}2.$
Abscissas of vertical jumps correspond to even values of~$d.$ Recall that\x{ the present} all our results apply to odd~$d$ only.

%----------------------------------------------------

It should be noted that whenever $n$ is large,
$\hmind/n$ approaches $1/4$ very quickly as $d$ increases starting from $d=1,$ as illustrated by Fig.~\ref{f:ThreeDom}b.

%-------------------------------------------------------------------------------------
\section{Security \lowercase{Gaps}}
\label{s:UncertainGap}

Since the regions of\x{ uniformly approving and uniformly disapproving configuration classes} uniform approval and uniform disapproval are symmetric with respect to the line $h=\frac n2,$ there exists a value $r$ that can be thought of as a \emph{$2$-way security gap\/}, in the sense that whenever the number of happy vertices deviates upwards (respectively, downwards) from $n/2$ by at least $r$, the resulting configuration class is guaranteed to be uniformly approving (resp., disapproving). In this section, we give a refined definition of this value and characterize the classes $\calGnd$ of regular graphs whose $2$-way security gap is sufficiently small.

In view of Theorem~\ref{t:Classes},
the value
\eq{e:UncertR}{
\tilde r(n,d) ~=~
\frac n2-\x{\beta_1(n,d)}\hmind ~=~
\hmand-\frac n2
}
can be called the \emph{uncertainty radius\/} of the class $\calGnd,$ since the inequality
$\hmind\le h\le\hmand$ that describes mixed classes $\calGndh$ can be written equivalently as
\eq{e:Uncert2a}{
\frac n2-\tilde r(n,d) ~\le~ h ~\le~ \frac n2+\tilde r(n,d)\,.\nonumber}%-

However, $\tilde r(n,d)$ can be negative: indeed, if $d=1$ or $d=n$ or $(n,d)=(5,3),$ then $\tilde r(n,d)=-\frac12,$ which corresponds to the case where $\calGnd$ contains no mixed classes $\calGndh.$ To fix this flaw, let us introduce
\eq{e:r}{r(n,d)~=~\tilde r(n,d)+\frac12 ~=~ \frac{n+1}2-\x{\beta_1(n,d)}\hmind,}
referred to as the \emph{$2$-way security gap}
of $\calGnd.$ It follows from Theorem~\ref{t:Classes} that
\eqss{
\label{e:Approv}%+
\calGndh \mbox{\;\:is a uniformly     approving class iff \;\:} h\!\!&>&\!\!\left\lfloor \frac n2+r(n,d)\right\rfloor;\\
\label{e:Disapprov}%+
\calGndh \mbox{\;\:is a uniformly  disapproving class iff \;\:} h\!\!&<&\!\!\left\lceil \frac n2-r(n,d)\right\rceil.
}

It follows from \eqref{e:UncertR} and \eqref{e:r} that
\eq{e:d-to-r}{r(n,d) ~=~ \left\lfloor\min\left\{\frac{(n+1)(d-1)}{4d};\;\frac{n-d}2\right\}\right\rfloor
\mathstrut-\indi{(5,3),(9,5)}{(n,d)},}
%+
where the minimum is equal to the first expression for low-degree graphs and to the second one for high-degree graphs.
Since by \eqref{e:Intersect} the minimum in $d$ of the lower bound of $\hmind$
is $(\check d(n),\check h(n))=(\frac{n+1}2,\frac{n+3}4),$ it holds that
\eq{e:delta}{0 ~\le~ r(n,d) ~\le~ \frac{n+1}2-\frac{n+3}4 ~=~ \frac{n-1}4\,.}%+

Given $n,$ it can be of interest
to find $d$ guaranteeing that the $2$-way security gap does not exceed a chosen threshold $r_0$:
\eq{e:SemiGap}{r(n,d) ~\le~ r_0\,.}%+

We say that $\calGnd$ is a \emph{class of\/ $2$-way $r_0$-\secure\ graphs} if it satisfies~\eqref{e:SemiGap}.
The following proposition is a direct consequence of \eqref{e:d-to-r} and the relevant definitions.

\begin{proposition}\label{p:Round}
Let
$\,0\le r_0\le\frac{n-1}4.$ Then $\calGnd$ is a class of\/ $2$-way $r_0$-\secure\ graphs if and only if
\eq{e:coR}{
\left[1\le d\le\frac{n+1}{n+1-4\left(r_0\mathstrut+\indi{(5,3),(9,5)}{(n,d)}\right)}
\mbox{\;\:~~~or~~~\,\;\:} n-2\left(r_0\mathstrut+\indi{(5,3),(9,5)}{(n,d)}\right)\le d\le n\right].\nonumber}%-
\end{proposition}

%-------------------------------------------------------------------------------------
\section{Vertex \lowercase{Degrees Allowing approval with Lowest Support}}
\label{s:MaxGap}

Let us fix the number of vertices $n$ and find the range $\calD_n$ of degrees $d$ for which there exist approving configurations with the minimum number of happy vertices,
\eq{e:hmin}{\hminn~=~\min_{1\le d\le n}\hmind\,.\nonumber}%-
 By \eqref{e:minh}, $\hminn\ge\frac{n+3}4.$
It follows from the symmetry of\x{ the bounds $\beta_1(n,d)$ and $\beta_2(n,d)$} $\hmind$ and $\hmand$ w.r.t. $h=\frac n2$ that the values of $d$ belonging to the range $\calD_n$ \emph{maximize\/} the requirements for $h$ to \emph{guarantee\/} approval.

To estimate $\hmind-\hminn,$ it is convenient to rewrite Proposition~\ref{p:Round} in a different form.
Due to \eqref{e:delta}, representation
\eq{e:delta1}{r_0 ~=~ \frac{n-1}4-\delta,\;\;\mbox{where}\;\;0\le\delta\le\frac{n-1}4}%+
is valid.
Using \eqref{e:delta1} and recalling that $\check h(n)=\frac{n+3}4$ according to \eqref{e:Intersect}, we obtain

\begin{corollary}\label{co:Round1}
Suppose that $\,0\le\delta\le\frac{n-1}4$ and $n\not\in\{5,9\}$. Then
\eq{e:coR1}{
\hmind-\check h(n) ~\le~ \delta
\mbox{\quad if and only if\;\quad} \frac{n+1}{2+4\delta} ~\le~ d ~\le~ \min\left\{\frac{n+1+4\delta}2;\;n\right\}\,.\nonumber}%-
Moreover$,$ all configuration classes ${\calG}_{5|d|h}$ are uniform\/$:$ disapproving when $h\le2$ and approving when ${h\ge3.}$
Classes of the form ${\calG}_{9|d|h}$ are mixed if and only if $d\in\{3,5,7\}$ and $h\in\{4,5\}.$ \x{ with $d\in\{1,9\}$ is uniformly disapproving $($uniformly approving$)$ when $h\le4$ $(h\ge5).$ If $d\in\{3,5,7\},$ then this class is uniformly disapproving $($uniformly approving$)$ when $h\le3$ $(h\ge6).$}
\end{corollary}

\begin{proof}
$\calGnd$ is a class of $2$-way $r_0$-\secure\ graphs whenever \eqref{e:SemiGap} holds, which by \eqref{e:r} amounts to $\frac{n+1}2-
\hmind\le r_0.$
Substituting \eqref{e:delta1} yields
$\hmind-\frac{n+3}4\ge\delta,$ which nonstrictly complements the premise $\hmind-\check h(n) \le\delta$ of Corollary~\ref{co:Round1}. Applying the same complement to the conclusion of Proposition~\ref{p:Round} and substituting \eqref{e:delta1} we obtain the conclusion of Corollary~\ref{co:Round1}. The statements on the classes ${\calG}_{5|d|h}$ and ${\calG}_{9|d|h}$ also follow from Proposition~\ref{p:Round} or Theorem~\ref{t:Classes}.
\end{proof}

We now find the range $\calD_n$ of degrees $d$ that, given $n,$ allow approval with the lowest possible support $\hminn$.
The results for different ``types'' of $n$ in terms of residuals are collected in Proposition~\ref{p:Foot}.

\begin{proposition}\label{p:Foot}
Let an odd $n$ be the order of a regular graph with loops and $n\not\in\{5,9\}$. Then $\hminn$ and the segments of $d$ satisfying conditions $\hmind=\hminn$ or\/ $\hmind\le \hminn+1$ are as given in~Table~$\ref{t:4m3}$.
%\vskip-.8em
\begin{table}[ht]
\begin{center}
{\footnotesize
\begin{tabular}{l|ccc|cc|cc}
$n$ expressed in terms of $\q\in\N$       &&$8\q-3,\,8\q-7\lK$\!\!\!     &&$8\q-1\lL$&&&$\lLl8\q-5$\!\!\!\!\!\!\!\!\!\\
\hline
$\hminn$                            &&$\frac{n+3}4\lK$\!\!\!&&\!\!\!$\frac{n+5}4\lL$&&&$\lLl\lB{\frac{n+\lU{5}}4}$\!\!\!\!\!\!\!\!\!\!\\
{$\calD_n\!=\!\{d \mid \hmind\!=\!\hminn\}$}\!\!\!     &&$\lD{\frac{n+1}2}\lK$\!\!\!&&\!\!\!$\left[\frac{n+5}4,\frac{n+3}2\right]\lL$&&&$\lLl\left[\frac{n+1}4,\frac{n+3}2\right]$\!\!\!\!\!\!\!\!\!\!\\
The length of the \x{$d$-}segment {$\calD_n$}\!\!\!    &&$0\lK$\!\!\!               &&\!\!\!$\frac{n+1}4\lL$                         &&&$\lLl\lD{\frac{n+5}4}$\!\!\!\!\!\!\!\!\!\!\\
\hline\hline
$n$  expressed in terms of $l\in\N$
&$12l-3\lK$                               &$12l-7\lK$                               &$12l-11$&$16l-1\lK$&$16l-9\lI$&$16l-5\lI$&$\lU{16l-13}$\\
\hline
$d:\hmind\le \hminn+1$
&\!\!\!$\left[\frac{n+9}6,\frac{n+5}2\right]$\lK & $\lU{\left[\frac{n+1}6,\frac{n+5}2\right]$\lK & $\left[\frac{n+5}6,\sta{\frac{n+5}2}\right]}$\!\!\!
&\!\!\!$\left[\frac{n+9}8,\frac{n+7}2\right]$\lJ & $\left[\frac{n+1}8,\frac{n+7}2\right]$\!\!\!&\!\!\!$\left[\frac{n+13}8,\frac{n+7}2\right]$\lJ& $\left[\frac{n+5}8,\sta{\frac{n+7}2}\right]$\!\!\!\\
The length of the $d$-segment
&$\frac{n+3}3\lK$                         &$\lU{\frac{n+7}3}\lK$                          &$\staa{\frac{n+5}3}$
&$\frac{3n+19}8\lJ$&$\frac{3n+27}8$&$\lD{\frac{3n+15}8}\lJ$&$\staa{\frac{3n+23}8}$\\
\hline
\end{tabular}
}
\caption{\emph{The segments of $d$ that\/$,$ given $n,$ allow approval with the lowest possible support~$h.$
In this table\/$,$ $x^*\stackrel{\mathrm{def}}{=}\min\{x;\,n\};$ $x^{**}\stackrel{\mathrm{def}}{=}\min\{x;\,n-1\}.$}
\label{t:4m3}}
\end{center}
\end{table}
\end{proposition}

\begin{figure}[ht]
\centerline{\small (a)\hspace{25em}(b)}
\centerline{
\includegraphics[height=15.0em]{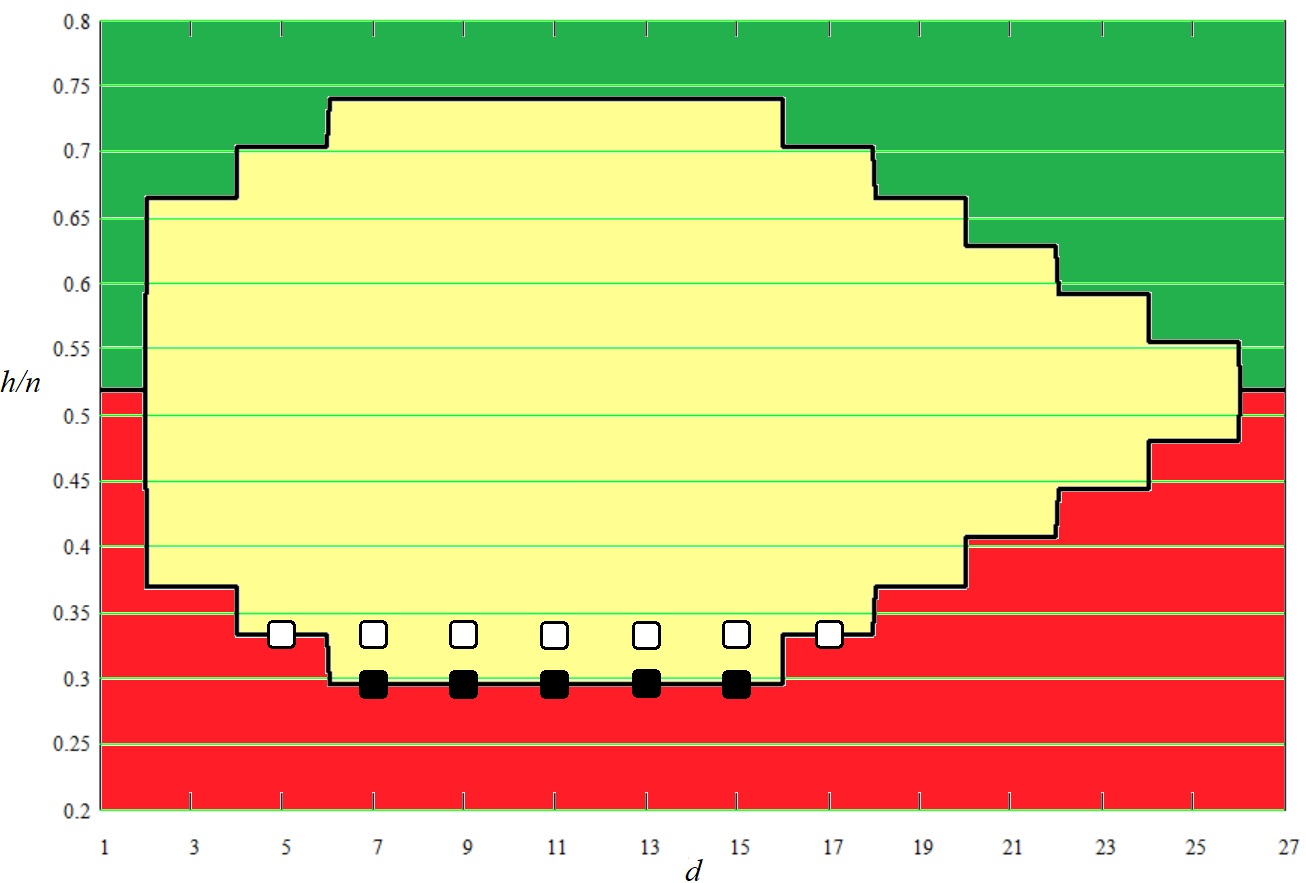}\;\; %
\includegraphics[height=15.0em]{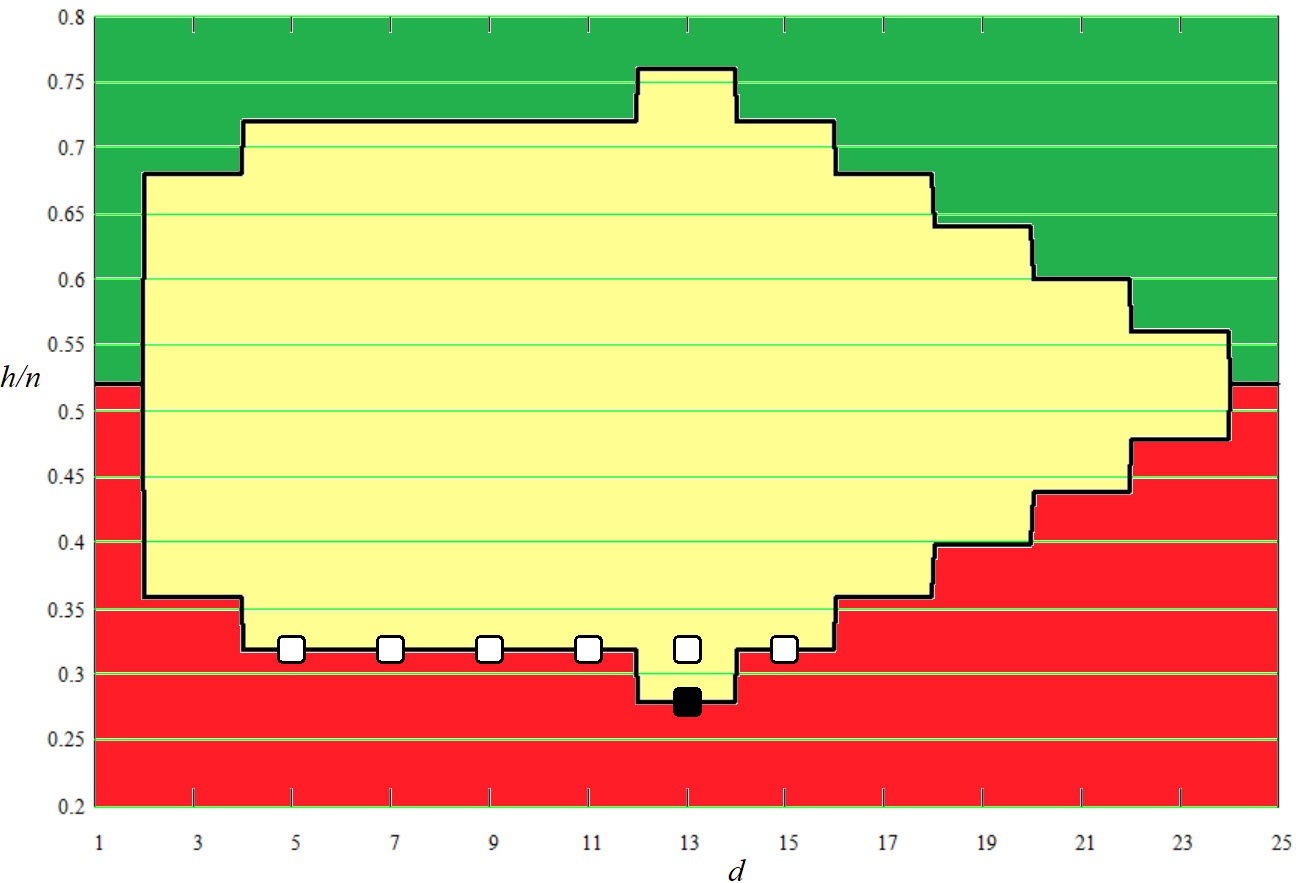} %
}
\caption{Regions of mixed (yellow), uniformly  disapproving (red), and uniformly  approving (green) classes $\calGndh$ of binary labeled $d$-regular graphs of order $n$ with $h$ happy vertices in the coordinates $d$ and $h/n$ for (a)~$n=27$ and (b)~$n=25.$ Classes $\calGndh$ that allow approval with minimum (respectively, minimum+1) support are shown by black (resp., white) squares.
\label{f:ThreeDom1}}
\end{figure}

%\vskip-1.4em
The cases of $n=16\!\cdot\!3-1=47,$ $n=8\!\cdot\!88-3=12\!\cdot\!59-7=701$ and $n=16\!\cdot\!2-5=27,$ $n=8\cdot4-7=12\cdot3-11=25$ are illustrated by Fig.~\ref{f:ThreeDom}a,b and Fig.~\ref{f:ThreeDom1}a,b, respectively.
Classes $\calGndh$ that allow approval with minimum (respectively, minimum+1) support are shown in Fig.~\ref{f:ThreeDom1} by black (resp., white) squares.

%----------------------------------------------------

\begin{proof}
Suppose that $n\not\in\{5,9\}$. We consider two cases.

1.~First consider the case of $n=4m-1,$ where $m\in\N.$ Then $\check h(n)=\frac{n+3}4=m+\frac12$ is half-integer. \x{By Theorem~\ref{t:Classes} it holds that
\eq{e:hmin1}{\hmind ~=~ \lceil\beta_1(n,d)\rceil.}}%+
We search $d$ such that $\hmind=\lceil\check h(n)\rceil=\frac{n+5}4.$ By
\eqref{e:hmimd} and Corollary~\ref{co:Round1}, for a half-integer $\check h(n)$ the following series of equivalences is true:
$$
\hmind=\lceil\check h(n)\rceil
~\Leftrightarrow~
\hmind-\check h(n) \le \frac12            ~\Leftrightarrow~
\frac{n+1}4\le d \le \frac{n+3}2           ~\Leftrightarrow~
m\le d \le 2m+1.
$$
If $m$ is odd, i.e.,\footnote{This $\q$ is not to be confused with the parameter of $\qq$-subdomination.} $n=8\q-5,\;\q\in\N,$ then the actual range of such odd $d$ is
$[m,\,2m+1]=\left[\frac{n+1}4,\,\frac{n+3}2\right],$
and the length of this $d$-segment is $\frac{n+5}4.$

Otherwise, if $m$ is even, i.e., $n=8\q-1,\;\q\in\N,$ then the actual segment of odd $d$ is
$[m+1, 2m+1]=\left[\frac{n+5}4,\frac{n+3}2\right]$,
whose length is $\frac{n+1}4.$

2. In the opposite case, where $n=4m-3\;(m\in\N),$ $\check d(n)=\frac{n+1}2=2m-1$ is odd and $\check h(n)=\frac{n+3}4=m=\hmin(n,\check d(n))$ is integer. Hence the class $\calG_{n|\check d(n)|\check h(n)}$ is nonempty; since $n\not\in\{5,\,9\},$ it contains an approving configuration by Theorem~\ref{t:Classes}.

If $d\ne\check d(n),$ then $\hmind>\check h(n)=m$ and therefore,
$$
\hmind\x{ ~=~ \lceil\beta_1(n,d)\rceil}~>~\lceil\check h(n)\rceil ~=~ \hmin(n,\check d(n)) ~=~ \hminn.
$$
Thus, $\check d(n)=\frac{n+1}2$ is the unique $d$ such that $\hmind=\hminn,$ which is illustrated by Fig.~\ref{f:ThreeDom1}b, where $n=25$.

3. Now let us find the range of $d$ that allow approval with a (low) support no more than $\hminn+1$ in the case of $n=4m-3.$
This amounts to $\hmind-\check h(n)\le 1,$ which by Corollary~\ref{co:Round1} is equivalent to $d\in\left[\frac{n+1}6,\min\!\left\{\frac{n+5}2;\,n\right\}\right].$
The right end of this segment, either  $\frac{n+5}2=2m+1$ or $n,$ is odd. To find the exact (odd)\x{ lower bound} left end, we have to consider three cases: $m=3l,$ $m=3l-1,$ and $m=3l-2$ ($l\in\N$). This is done straightforwardly; the results are presented in Table~\ref{t:4m3}. The case of $n=4m-1$ is considered similarly by using Corollary~\ref{co:Round1} with $\delta=1.5.$
\end{proof}

It follows from Proposition~\ref{p:Foot} that
the length of the $d$-segment on which approval can be implemented with the least possible support is $0$ when $\modd{n}{4}=1$ (except for $n\in\{5,9\}$), $\frac n4+0.25$ when ${\modd{n}{8}=7},$ or $\frac n4+1.25$ when $\modd{n}{8}=3.$
The length of the segment where $d$ satisfies condition $\hmind\le \hminn+1$ varies from $\frac n3+1$ to $\frac{3n}8+3.375.$ The left boundary of this segment, depending on $n,$ slightly exceeds $\frac n6$ or even $\frac n8.$
Sets of $d$ corresponding to other conditions of this kind can be found using Corollary~\ref{co:Round1} or Proposition~\ref{p:Round} similarly.

\begin{figure}[ht]
\centerline{\small (a)\hspace{24em}(b)}
\centerline{
\includegraphics[height=15.0em]{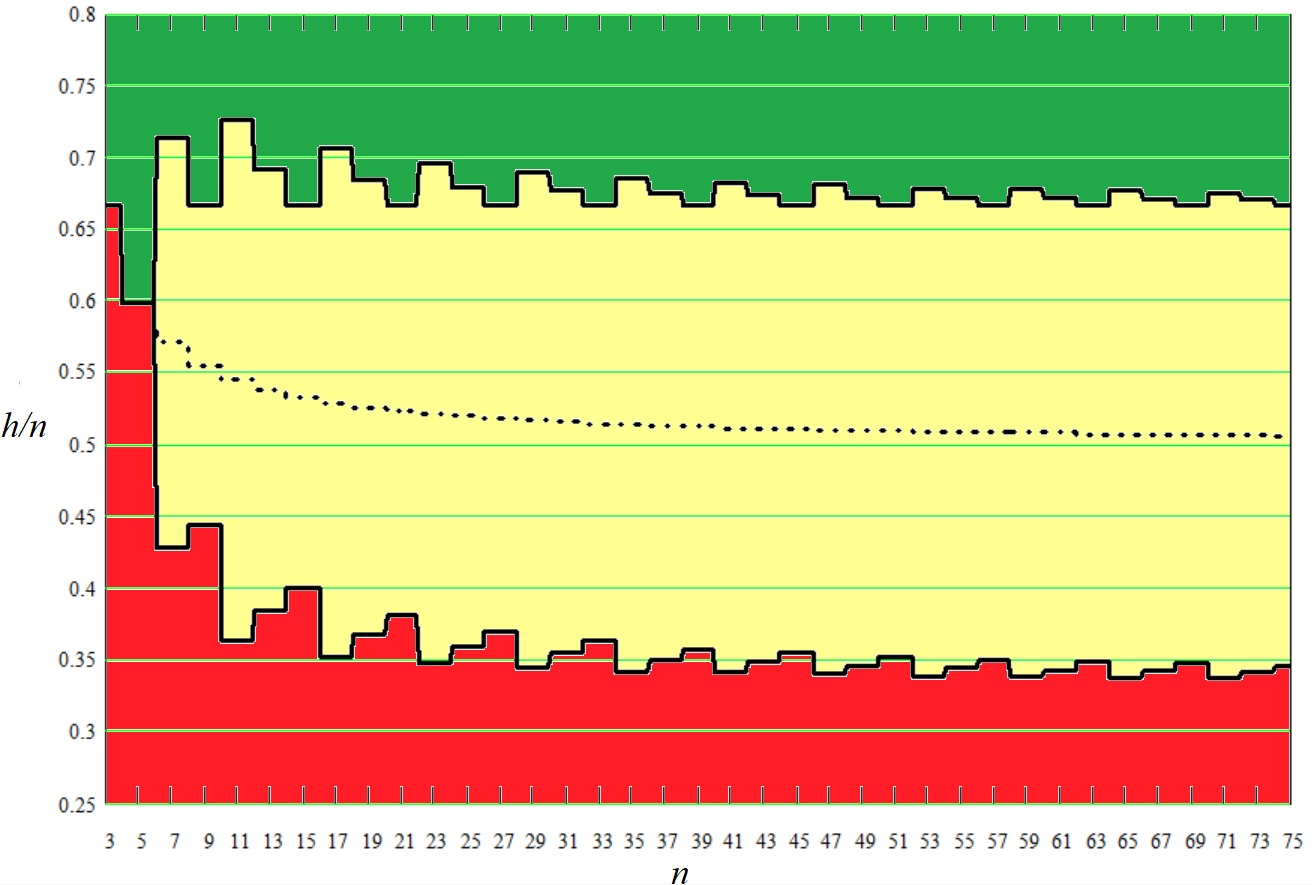}\;\; %
\includegraphics[height=15.0em]{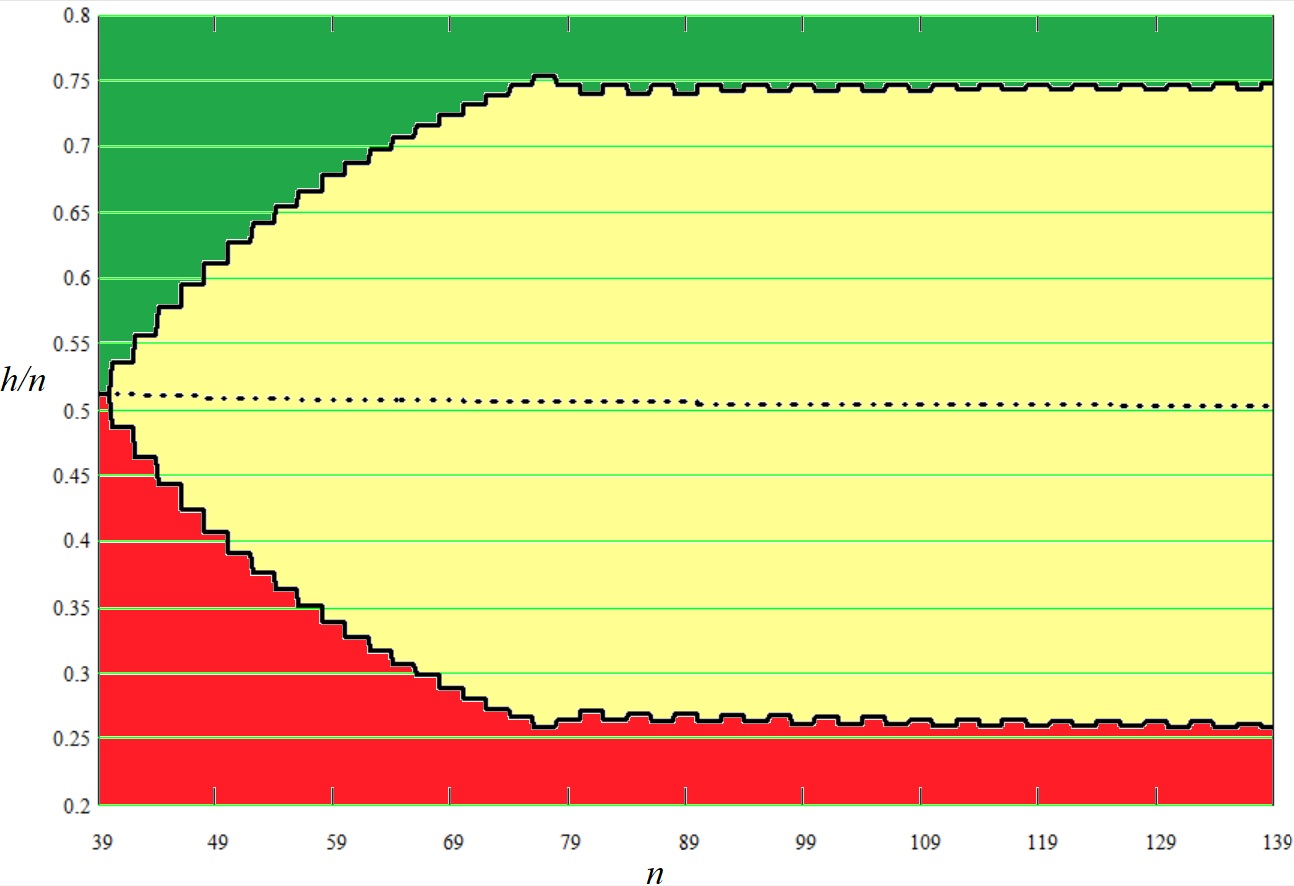} %
}
\caption{Regions of mixed (yellow), uniformly  disapproving (red), and uniformly  approving (green) classes $\calGndh$ of binary labeled $d$-regular graphs of order $n$ with $h$ happy vertices in the coordinates $n$ and $h/n$ for (a)~$d=3$ and (b)~$d=39$. The dependence of $h/n=\frac {n+1}{2n}$ on $n$ is shown by a dotted curve.\label{f:ThreeDom2}}
\end{figure}

The\x{ same} regions of mixed, uniformly  disapproving, and uniformly  approving classes $\calGndh$\x{ of binary labeled $d$-regular graphs of order $n$}
are shown in Fig.~\ref{f:ThreeDom2} in the coordinates $n$ and $h/n$ for fixed $d=3$ and $d=39.$
The points on the separating boundaries are symmetric w.r.t. $h/n=\frac {n+1}{2n}$; the dependence of this value on $n$ is shown by a dotted curve.

%----------------------------------------------------

The left part of diagram (b) corresponds to high-degree graphs with $n$ satisfying $39=d\le n\le2d-1=77.$ For them, the boundary separating uniformly  disapproving and mixing classes decreases hyperbolically with increasing $n$ due to the local support inequality~\eqref{e:LocalIneq}: $\frac hn\ge\frac{d+1}{2n}.$ The right part refers to low-degree graphs with $n>77$ for which this boundary decreases hyperbolically too, but according to the global support inequality~\eqref{e:GlobalIneq1}, $\frac hn\ge\frac14\big(1+n^{-1}\big)\big(1+d^{-1}\big)$, converging to $\frac14\big(1+39^{-1}\big)=\frac{10}{39}\approx 0.2564$ as $n\to\infty.$

In the left diagram (a), the difference between these parts is not noticeable due to corrections related to the integer values of $h,$ odd values of $n,$ and the exceptions $(n,d,h)\in\{(5,3,2),(5,3,3)\}$
of Theorem~\ref{t:Classes}. Thus, with $d=3$ we have $h\ge2$ for $n=3,$ $h\ge3$ for $n=5,$ and $h\ge\frac{n+1}3$ for low-degree ($n>5$) graphs; $\lim_{n\to\infty}\frac hn=\frac13$.

\begin{figure}[ht]
\centerline{\small (a)\hspace{24em}(b)}
\centerline{
\includegraphics[height=15.9em]{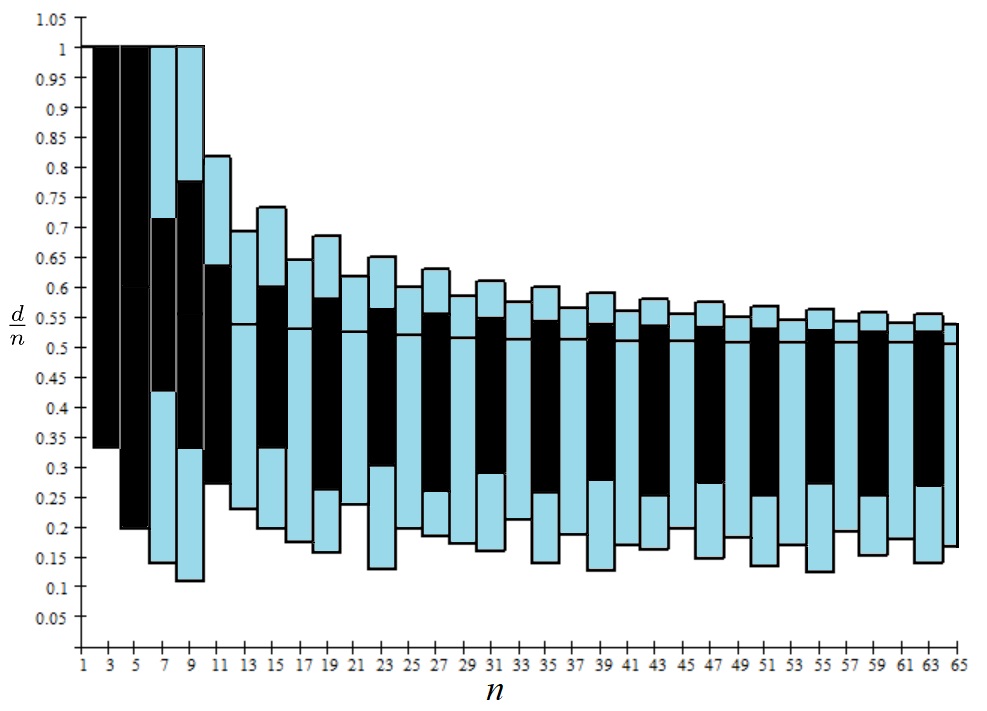}\;\; %
\includegraphics[height=15.9em]{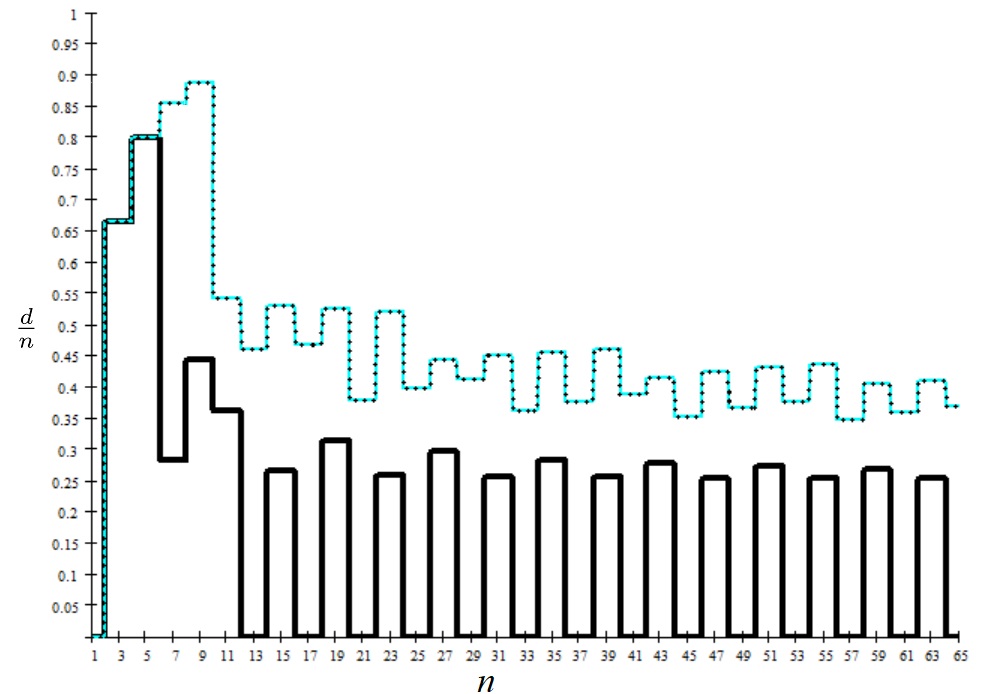} %
}
\caption{(a) Segments of the degree $d$ values such that $\hmind=\hminn$ (foreground, in black), and segments of $d$ values with $\hmind\le \hminn+1$ (background, in cyan), measured relative to $n$ and plotted against~$n;$ (b)~The length of these segments\x{ of $d$ such that $\hmind=\hminn$ (black) and $\hmind\le \hminn+1$} divided by $n$ (in solid black and dotted cyan, respectively).}
\label{f:segments}
\end{figure}

\smallskip
Figure~\ref{f:segments} directly shows the segments of the values of $d$ that enable approval with the minimum or `minimum+1' support along with the relative length of these segments.

Proposition~\ref{p:Foot} and Corollary~\ref{co:Round1} imply the following.

\begin{corollary}\label{co:nmin}
For any odd\/ $n\in\N,$
$\hminn=\frac14(n+2+\modd{n}{4})+\indi{5,9}{n}\x{;\;\hmin(5)=3;\;\hmin(9)=4}.$
\end{corollary}

The following corollary gives general expressions for the values presented in Table~\ref{t:4m3} and related to the $d$-segments of approval with low support.

\begin{corollary}\label{co:segm}
For $m\ge0,$ let $L_m(n),R_m(n),$ and $D_m(n)$ denote the left end\/$,$ the right end\/$,$ and the length of the segment of the $d$ values such that
$\hmind\le \hminn+m.$ Then for any odd\/ $n\in\N\!:$
\eqss{\nonumber
L_0(n)&=&\frac18\big((\modD{n}{4}-1)(\modD{n}{8}-n-4)+4n+4\big)-2\!\cdot\!\indi{5,9}{n}\,;\\\nonumber
R_0(n)&=&\frac12(n+\modD{n}{4})+2\!\cdot\!\indi{5,9}{n}\,;\\\nonumber
D_0(n)&=&\frac18(\modD{n}{4}-1)(n+8-\modD{n}{8})+4\!\cdot\!\indi{5,9}{n}\,;\\\nonumber
L_1(n)&=&\frac1{16}(\modD{n}{4}-1)(n+\modD{n}{16}-2\modD{n}{8}+8)\\\nonumber
&&~~ -\frac1{12}(\modD{n}{4}-3)(n-4\modD{n}{3}+9)-2\!\cdot\!\indi{9}{n}\,;\\\nonumber
R_1(n)&=&\min\left\{\frac12(n+\modD{n}{4}+4);\;n\right\}+2\!\cdot\!\indi{9}{n}\,;\\\nonumber
}
\eqss{\nonumber
D_1(n)&=&\min\bigg\{\frac1{16}(\modD{n}{4}-1)(3n+2\modD{n}{8}-\modD{n}{16}+20)\\\nonumber
&&~~~~~ -\frac16(\modD{n}{4}-3)(n+2\modD{n}{3}+3);\;n-1\bigg\}+4\!\cdot\!\indi{9}{n}\,.
}
\end{corollary}

These\x{ results} expressions are derived directly from Proposition~\ref{p:Foot} and Corollary~\ref{co:Round1}\x{ straightforwardly}.

%-------------------------------------------------------------------------------------
\section{Alternative formulation: Majority \lowercase{Domination on Regular Graphs}}
\label{s:MajorityRegular}

The question regarding the minimum value of $h$ such that $\calGndh$ contains approving configurations can be reformulated in terms of majority domination.
The first survey on domination in graphs \cite{CockayneHedetniemi77} appeared in 1977 and contained 20 references.
A \emph{majority $[$dominating\/$]$ function\/} \cite{Broere95majority} on $V=V(G)$ is an opinion function $f: V\to\{1,-1\}$ such that $f(N_u)\ge1$ for \emph{at least half\/} of the vertices $u$ in~$V$.
A \emph{strict majority function\/} is an opinion function $f: V\to\{1,-1\}$ such that $f(N_u)\ge1$ for \emph{more than half\/} of the vertices $u$ in~$V$.
For a positive integer $\qq,$ a \emph{$[$signed\/$]$ $\qq$-subdominating function\/} of $G$ is an opinion function $f: V\to\{-1, 1\}$ such that $f(N_u)\ge1$ for \emph{at least $\qq$} vertices $u$ of~$G.$ The latter concept reduces to the previous two when $\qq=\lceil\card{V}/2\rceil$ and  $\qq=\lceil(\card{V}+1)/2\rceil,$ respectively.

The [weak] \emph{majority domination number\/}         $\maj(G)$           \cite{Broere95majority},
the \emph{strict majority domination number\/}         $\smaj(G)$          \cite{HenningHind98strict,KangShan00lower}, and
the \emph{$[$signed\/$]$ $\qq$-subdomination number\/} $\gamma^{}_{\qq s}(G)$ \cite{Cockayne1996}
of $G$ are the minimum possible weights of
a majority function,
a strict majority function, and
a [signed] $\qq$-subdominating function
on $V(G),$ respectively.

N.~Alon proved (see~\cite{Broere95majority,Cockayne1996}) that the majority domination number $\maj(G)$ of a connected graph $G$ is at most~2. Moreover, $\maj(G)$ does not exceed $1$ when the order $n$ of $G$ is odd~\cite{Cockayne1996,HenningHind98strict}.
Obviously, in the latter case,
$\maj(G)=\smaj(G).$
The decision problem corresponding to computing $\maj(G)$ is NP-complete~\cite{Broere95majority}.

Let
\eq{e:maj(nd)}
{\maj(n,d)~=~\min_{G\in\calGnd}\maj(G);\quad \maj(n)~=~\min_{1\le d\le n}\maj(n,d)\,,}
{i.e., $\maj(n,d)$ is the minimum $f(V)$ over majority functions $f$ on~$V(G)$ for $G\in\calGnd.$}
Observe that
\eq{e:f(V)}{f(V)~=~h-(n-h)~=~2h-n\,,}%+
where $h$ is the number of happy vertices.
Hence, for an odd\x{order } $n,$ it holds that
\eq{e:majmin}{\maj(n,d)~=~2\hmind-n\,.}%+
Indeed, in the case of odd $n$, any weak or strict majority function $f$ ensures that $f(N_u)\ge1$ for more than half of the vertices $u$ of~$G,$ i.e., that more than half of the vertices are proponents, which means that the graph labeled with $f$ is an approving configuration. Therefore, finding $\maj(n,d),$  $\smaj(n,d)$ (defined in the same way as in~\eqref{e:maj(nd)}), and $\hmind$ involves minimization over the same sets of graphs and opinion functions~$f.$

Henning~\cite{Henning96domination} and independently Holm~\cite{Holm01majority}\x{ Received 8 October 1996}
proved the following inequality (presented here in our notation): for a $d$-regular graph of order $n$ with loops and odd $d\ge3,$
\eq{e:HenHol}{\maj(G) ~\ge~ -\frac n2\cdot\frac{d-1}d~.}%+
The same lower bound can be obtained from a more general result \cite{KangShan00lower} on graphs with specified minimum and maximum vertex degrees.

On the other hand, for odd $n,$ it follows from \eqref{e:GlobalIneq} that
\eq{e:g_maj}{\maj(G)~=~2\hmind-n~\ge~ 2\frac{n+1}2 \cdot \frac{d+1}{2d}-n ~=~ -\frac{n(d-1)}{2d}+\frac12+\frac1{2d}~,}%+
which is stronger in some cases than~\eqref{e:HenHol}, although it is claimed \cite{Ungerer96PhD,Henning96domination,HattinghUngererHenning98,Harris03PhD,KangQiaoShanDu03} that the bound~\eqref{e:HenHol} is sharp or best possible.
Technically,
it is sharp in a weak sense, that is,
it cannot be improved for \emph{some\/} values of $n$ and~$d.$ Proposition~\ref{p:Nonsharp} below states that for any odd $d>1$ there are infinitely many $n$ such that
bound \eqref{e:HenHol} can be improved for the class $\calGnd$ of $d$-regular graphs with loops on $n$ vertices. Moreover, the inaccuracy of \eqref{e:HenHol} reaches $\frac{n+1}2$ in the case of high-degree graph. The bounds obtained in this paper are sharp for any pair $(n,d)$ with odd~$n.$

\begin{proposition}\label{p:Nonsharp}
There is no odd $d\ge3$ such that bound \eqref{e:HenHol} is sharp for order $n$ $d$-regular graphs with loops and all odd $n\ge d$.
The inaccuracy of \eqref{e:HenHol} reaches $1$ for low-degree graphs and $\frac{n+1}2$ for high-degree graphs.
\end{proposition}

\begin{proof}
For any odd $d\ge3$ let $n=2dl+3,\,l\in\N.$ Consider
the class $\calGnd$ of $d$-regular graphs with loops on $n$ vertices, which are low-degree. For this class, \eqref{e:g_maj} takes the form
\eq{e:h-glo}{
\maj(n,d) ~\ge~ -\frac{nd-n-d-1}{2d} ~=~ -(dl-l+1)+\frac2d~,
}%+
while \eqref{e:HenHol} provides
\eq{e:h-hen}{\maj(n,d) ~\ge~ -(dl-l+1)-\frac12\left(1-\frac3d\right).}%+

It follows from \eqref{e:h-glo} and \eqref{e:h-hen} that the value $\maj(G)=-(dl-l+1)$ is prohibited by the global support inequality \eqref{e:GlobalIneq} for all $G\in\calGnd,$ but is allowed by~\eqref{e:HenHol}. Hence, bound \eqref{e:HenHol} is not sharp for the pairs $(n,d)$ under consideration; as $d\ge3,$ its inaccuracy is~$1.$

Now consider high-degree graphs, for which $n/2<d\le n$. The local support inequality \eqref{e:LocalIneq} is equivalent to
\eq{e:Locg1}{\maj(n,d)~=~2\hmind-n ~\ge~ d+1-n\,.}%+

The difference between the bounds \eqref{e:Locg1} and \eqref{e:HenHol} is
\eq{e:Diffe}{
(d+1-n)+\frac{n(d-1)}{2d} ~=~ (d+1)\left(1-\frac n{2d}\right) ~\le~ (n+1)\left(1-\frac12\right)~=~\frac{n+1}2\,,\nonumber}%-
where $\frac{n+1}2$ is reached when $d=n.$
\end{proof}

Note that the more accurate lower bound
$$
\maj(G) ~\ge~ \frac{n+1}2\cdot\frac{d+1}d-n\,,
$$
which coincides with \eqref{e:g_maj} and is equivalent to the global support inequality, can be derived from a result of \cite{Ungerer96PhD,HattinghUngererHenning98} (cf. \cite[Theorem~4.15]{Hattingh98majority}), which in turn follows from the results in \cite{ChangLiawYeh02,KangQiaoShanDu03} and Theorem~2 in \cite{ChenSong08lower} on $k$-subdomination numbers, by substituting $k=\frac{n+1}2$ (in the notation of\x{ \cite{Ungerer96PhD} and} \cite{HattinghUngererHenning98},\x{ $t=\frac{n+1}{2n}$ and} $q=\frac{n+1}{2n}$\x{ resp.}) in the case of odd~$n.$
As shown in Section~\ref{s:Constru},
this bound is sharp for low-degree graphs. For high-degree graphs it is not sharp, and its inaccuracy in bounding $\hmind$ can be estimated by $\left(1-\frac{n+1}{2d}\right)\!\frac{d+1}2,$ which reduces to $\frac{n^2-1}{4n}$ when $d=n$. This follows from considering the difference between the bounds \eqref{e:LocalIneq} and \eqref{e:GlobalIneq} and the propositions on the sharpness of bound \eqref{e:LocalIneq} proved in Section~\ref{s:Constru}.

The sharpness of the bound \eqref{e:g_maj} generalized to the $\qq$-subdominating functions is claimed in~\cite[Theorem~4.1]{Ungerer96PhD}. However, to prove this for a given $d,$ the author constructs a\x{\x{ specific} $(d-1)$-regular graph\x{ $G$} without loops} graph whose order $n$ is determined by the parameter $\qq$ of subdomination and is divisible by~$d$ (like in \cite{Henning96domination}) so this is sharpness in a weak sense.
Similarly, the sharpness of the bounds given in \cite[Theorem~2]{ChenSong08lower} is established only
for some values of the parameters, namely, when $\qq=n,$ or in other words, for $n$-subdomination (also called \emph{signed domination}). In Section~\ref{s:Constru}
we proved
that the proposed bounds
are attainable for all odd $n$ and $d\le n$.

A number of additional inequalities related to majority domination and $\qq$-subdomination\x{ are presented} can be found in~\cite{KangShan20signed}.
For other results applicable to regular graphs with some variations of the majority domination model we refer to
\cite{Dunbar96minus,HaynesHedetniemi98book,LiuSunTian02lower,Xing05signed,HarrisHattinghHenning06,KangShan07survey,HamidPrabhavathy16majority,Caro18effect}.

In terms of the weight of majority [dominating] functions, the core Problem~\ref{pm:1} considered in this paper can be reformulated as follows.

\begin{manual}[problem]{\ref{pm:1}\thmprimeprime}
\label{pm:1''}
Given odd $n,d\in\N$ such that
$d\le n,$ find\/$:$
\begin{description}
\item{$(1)$}
$\maj(n,d);$
\item{$(2)$}
all $l\in[-n,n]\cap\N$ such that $f(V)=l$ with $f\!:V(G)\to\{1,-1\}$ and $G\in\calGnd$ implies that $f$ is a majority function\/$;$
\item{$(3)$}
all $l\in[-n,n]\cap\N$ such that $f(V)=l$ with $f\!:V(G)\to\{1,-1\}$ and $G\in\calGnd$ implies that $f$ is not a majority function\/$.$
\end{description}
\end{manual}

\x{A solution to this problem is given by Corollary~\ref{co:Classes} in Subsection~\ref{s:Classes}.}
\smallskip

The solution to Problem~\ref{pm:1}, given by Theorem~\ref{t:Classes},
is now reformulated in term of $f(V),$ the (odd) weight of the corresponding opinion function
related to $h$ by \eqref{e:f(V)} and $\maj(n,d)$.
Together with Eq.~\eqref{e:majmin}, this yields the following corollary,
which solves Problem~\ref{pm:1''}.

\begin{corollary}\label{co:Classes}
Suppose that $n,d\in\N$ are odd\/$,$ $d\le n,$ and $f(V)=\sum_{v\in V}f(v)$ is the weight of an opinion
function $f: V\to\{1,-1\},$ where $V=V(G)$ is the vertex set of a graph $G\in\calGnd.$
Then\/$:$

\smallskip\noindent
$\displaystyle
(1)~\maj(n,d)
\,=\,2\left(\left\lceil\frac{d+1}{2}\cdot\max\left\{\frac{n+1}{2d}\,,\,1\right\}\right\rceil+\indi{(5,3),(9,5)}{(n,d)}\!\right)-n\,;$

\smallskip\noindent
$(2)$~every odd number in
$(-\maj(n,d),\,n]$
is the weight of some majority function $f$ and is not the weight of any non-majority function~$f$$;$

\smallskip\noindent
$(3)$~every odd number in $[-n,\,\maj(n,d))$ is the weight of some non-majority function $f$ and is not the weight of any majority function $f$.
\end{corollary}

In other words, due to Corollary~\ref{co:Classes},
the region of mixed classes in terms of $f(V)$
is $[\maj(n,d),$ $-\hspace{.07em}\maj(n,d)],$ so that there are no mixed classes when $\maj(n,d)>0$ and the regions of uniformly approving and uniformly disapproving classes are symmetric w.r.t. the line $f(V)=0.$

Finally, we reformulate the results of Proposition~\ref{p:Foot} and Corollary~\ref{co:nmin} that involve $\hmind$ or $\hminn$ in terms of
(an explicit form of)
$\maj(n,d)$ and $\maj(n)$.
Using \eqref{e:majmin} we obtain the following.

\begin{corollary}\label{co:Foot}
For any odd\/ $n\in\N,$
$\maj(n)=\frac12(2+\modd{n}{4}-n)+2\!\cdot\!\indi{5,9}{n}\x{;\;\maj(5)=1;\;{\maj(9)=-1}}.$
\x{If $n=4m-1\,(m\in\N),$ then $\maj(n)=\frac{5-n}2;$
if $n=4m-3\,(m\in\N),$ then $\maj(n)=\frac{3-n}2.$}

Conditions $\big[d:\hmind=\hminn\big]$ and\/ $\big[d:\hmind\le \hminn+1\big]$ in Table~$\ref{t:4m3}$ are equivalent to $\big[d:\maj(n,d)=\maj(n)\big]$ and $\big[d:\maj(n,d)\le\maj(n)+2\big],$ respectively\/$,$\x{ amount to}  and can be replaced by them.
\end{corollary}

\begin{figure}[ht]
\centerline{\small (a)\hspace{24em}(b)}
\centerline{
\includegraphics[width=22.1em]{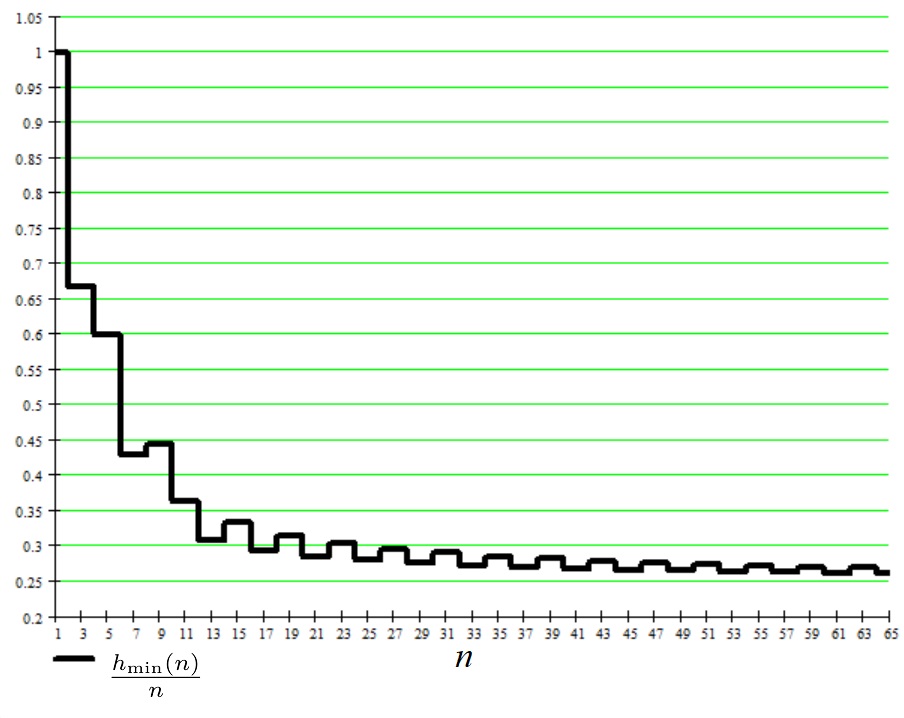}\;\; %
\includegraphics[width=22.1em]{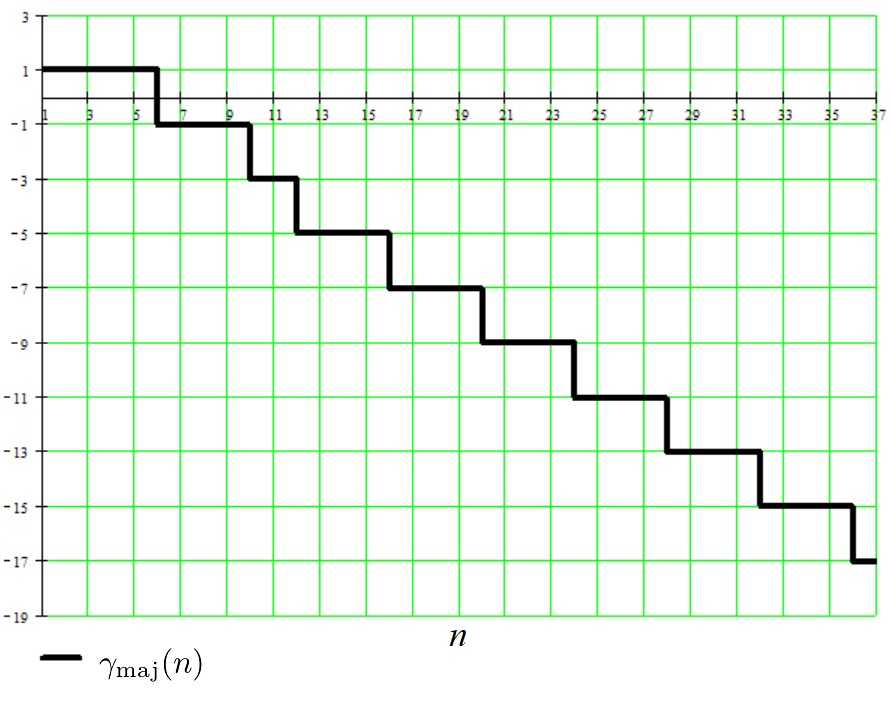} %
}
\caption{Dependence of $\frac{\hminn}n$\x{ (a)} and $\maj(n)=\smaj(n)$\x{ (b)} on $n\in\Nodd=\{1,3,5,\ldots\}.$} \label{f:hminGamma}
\end{figure}

Figure~\ref{f:hminGamma} shows how $\frac{\hminn}n$ and $\maj(n)$ depend on~$n\in\Nodd=\{1,3,5,\ldots\}.$ It can be noted that $\maj(4m+1)=\maj(4m+3)=\maj(4m+5)+1$ for any $m\in\N\!\smallsetminus\!\{1,2\}.$

%----------------------------------------------------
\section{Discussion}
\label{s:Discus}

%-------------------------------------------------------------------------------------
\subsection{Connections to \lowercase{Related Work}}
\label{ss:Related}
%-------------------------------------------------------------------------------------

Related works in the framework of graph domination have been discussed in Section~\ref{s:MajorityRegular}. Now we\x{ consider the} turn to studies carried out\x{ within} in different\x{ other, alternative} frameworks.

%-------------------------------------------------------------------------------------
\paragraph{The power of small coalitions.}

\begin{figure}[ht]
\centerline{
\includegraphics[height=5em]{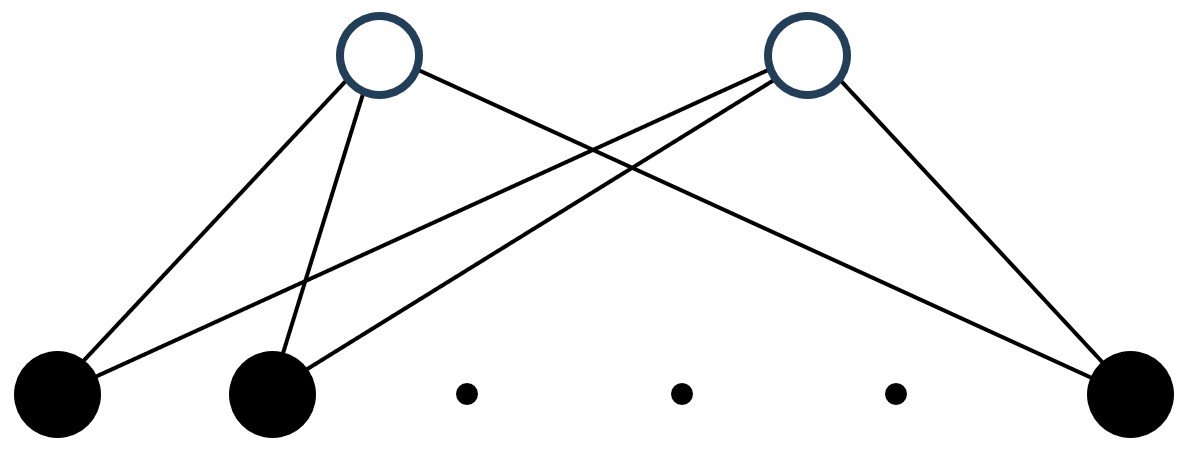}
}
\caption{A complete bipartite graph with two vertices that control the majority of local majorities \label{f:bipart}}
\end{figure}

The classical lower bound of $\frac14$ for the relative support of decisions made by a majority of majorities is valid in the special case where all local voting bodies have the same size and do not overlap. In the
more general case (of overlapping voting bodies of different sizes), the support can be arbitrarily low.
This is illustrated by Fig.~\ref{f:bipart} \cite{Peleg02local}, where the majority of majorities makes a decision supported by only two happy vertices (depicted in white) in the presence of arbitrarily many sad ones (shown in black). Here, every vertex is supposed to belong to its own neighborhood.

The same low support is sufficient for making decisions in the case of voting bodies \emph{of the same size\/} if \emph{non-symmetric connections\/} are allowed. Indeed, if each of the two white vertices in Fig.~\ref{f:bipart} takes into account the opinions of itself and exactly two black vertices, then all voting bodies have size 3 and the majority of majorities still adopts the proposal put to vote.

On the other hand, a negligible minority of $2\sqrt{n}$ happy voters can provide a $2:1$ majority in \emph{all\/} local voting bodies of \emph{varying size\/} with \emph{symmetric connections\/}~\cite{LinialPeleg93}.

In these examples, the members of small coalitions controlling the majority have much more outgoing influence arcs than the other voters.
This suggests that the influence of a vertex is determined by its out-degree and the influence of a coalition is determined by the sum of out-degrees of its members.

Combining the condition of \emph{symmetric connections\/} with the condition of \emph{equal-sized\/} local voting bodies\x{, we obtain} gives a rather strong requirement.\x{, so it is interesting to find out to what extent} In the present paper we studied the extent to which it limits the power of small coalitions.

A structure satisfying this requirement and assuming that each voter belongs to their own voting body is modeled by a regular graph with loops. The main goal of\x{ the current paper is} this paper was to study the relationship between the acceptance of proposals and the share of their support in such structures.

The properties of majority polls determine in turn the characteristics of the dynamical processes governed by majority~\cite{Peleg98size,Auletta15minority,Gartner18majority,Avin19majority}.

Applications of such studies are not limited to traditional fields such as political, social, and management sciences, they also arise in modern swarm robotics \cite{Scheidler15kUnanimity}, collaborative teleoperation \cite{Goldberg03CollabTeleoper}, internet of things \cite{Makhdoom19IoTBlockchain}, machine learning \cite{Pourdarbani20MajorityNN}, network decontamination \cite{Luccio07network}, database management algorithms \cite{Wang07majoritywatermarking}, and distributed computing  \cite{HendrickxOlshevskyTsitsiklis11,Cruciani21}, where aggregating the results of local or global majorities is a means of overcoming failures. Another important application is the study and control of social networks~\cite{Auletta15minority,Avin19majority}.

\paragraph{A social choice perspective.}
In terms of social choice theory, a faithful decision-making policy approves a proposal only if it is a \emph{Condorcet winner\/} over the \emph{status quo} in two-alternative voting. This simply means that a\x{ simple} majority of voters prefer the proposal to the current state of society.
It has been shown that a two-tier majority is much less stable under\x{ perturbations} random vote corruption than a simple majority (see, e.g., \cite{Heilman21Stable,MouzonLaurent20USElect,Kalai10Noise}).

The ability of two-tier majority to represent all eligible voters is minimal when the group size and the number of groups are equal to each other~\cite{KaoCouzin19modular,BottcherKernell22Condorcet}. However, in the situations where only one of two possible answers is correct, the opinions are statistically independent (like in the Condorcet Jury Theorem), and voter abstention is allowed and correlated with competency within groups, two-tier majority often outperforms direct voting~\cite{BottcherKernell22Condorcet}.

In \cite{FeixLepelleyMerlin21} various methods are discussed for assigning weights (number of votes, mandates, etc.) in the second\x{ tier} stage of a two-stage voting system to voting bodies of different sizes, taking into account the faithfulness criterion. One of the conclusions of the authors is that the square root rule (assigning weight $\sqrt{n_i}$ to a voting body of size~$n_i$), ``which stood for a long time as the only normative recommendation for voting in federations, can be seriously contested'' by some other methods with weights of the form~$n_i^\delta,$ where $\delta\ge0.$

%----------------------------------------------------
\paragraph{Two-stage majority in the ViSE model.}

In the ViSE (Voting in Stochastic Environment) model \cite{Che06ARC,MaksChe20}, each voter is characterized by their capital; proposals put to the vote are formulated in terms of voters' capital gains and generated stochastically. Voters can have different voting strategies, one of which is the Group Principle~A \cite{Che06ARC}, according to which all members of a group support a proposal if and only if it is beneficial for the majority of group members. The final decision is made by the majority of all voters, so in the case of equal-sized groups, it is equivalent to
the majority of majorities. Usually the group strategy is more profitable than the selfish one, but it does not protect voters from the systematic loss of capital in an unfavorable environment~\cite{CheMal18opt}. Furthermore, under certain circumstances, small groups can be more effective than large ones~\cite{CheLog10SCW}.

%----------------------------------------------------
\paragraph{Majority on circulant graphs.}

``Do local majorities force a global majority?'' by Fishburn, Hwang, and Lee \cite{Fishburn86local} appeared in 1986, almost ten years before~\cite{Broere95majority} and the close paper~\cite{DunbarHedeHenSla95} on majority domination in graphs. The question refers to the situation where all happy (in \cite{Fishburn86local}, ``white'') vertices must be proponents, while a vertex is a proponent if and only if there are at least $c>0$ more happy vertices than sad vertices in its neighborhood. This neighborhood does not contain the vertex itself, which corresponds to graphs without loops in our framework.
The structure of links is specific and generally corresponds to a directed graph, but in the symmetric case, it is representable by a simple circulant graph (the so-called Harary graph). Although this case is special, it is an important case, because many graphs that ensure approval with minimum support are just circulant and even Harary graphs (cf.\x{ Item~4 of}\ the proofs of Propositions~\ref{p:High1}--\ref{p:High3} and \ref{p:Low} above).

While in the graph domination literature the focus of attention is the value of $h-s,$ in \cite{Fishburn86local}, $h/s$ is estimated. The authors prove that for $d$-regular approving configurations (in which all happy vertices must be proponents) in the symmetric case, $h/s$ always exceeds~$1.$ Moreover, $\frac hs\ge\frac{d+c}{d-c}.$ They also obtain some upper bounds (later improved by Lison\v{e}k \cite{Lisonek95local}) for the minimum of~$h/s.$

A number of advances in the study of this problem were made by Woodall~\cite{Woodall92local}, including some extensions to strongly connected digraphs and connected graphs with minimum degree at least~$2.$

Despite the similarity of the problems studied, this research line is almost ignored in the graph domination literature.

%-------------------------------------------------------------------------------------
\subsection{Technical summary}
\label{ss:summary}
%-------------------------------------------------------------------------------------

In this paper, we studied two-stage majority approval on regular graphs with loops. \x{Suppose there is}Let there be a proposal whose implementation would make each vertex happy or sad. On the first stage of the decision making process, every vertex becomes\x{ is recognized as} a proponent or an opponent. A proponent is a vertex whose neighborhood contains more happy vertices than sad vertices; in the opposite case, it is an opponent. On the second stage, the proposal is accepted if and only if a majority of vertices are proponents.

To eliminate symmetric situations, where exactly half of the vertices are proponents, we considered the case, where the order $n$ of the graph is odd.
In the presence of loops, this implies that the degree $d$ of each vertex is odd as well: otherwise, one of the conditions of the Lemma~\ref{lem:reg-graph}
is violated. This means in turn that a vertex neighborhood cannot split equally into factions of happy and sad vertices\x{ too}.

The focus of the study is the relationship between $n$, $d,$ and the number of happy vertices $h$ that cause approval. The whole range of $0\le h\le n$ splits into three parts: high $h$ values guarantee approval for any $d$-regular graph of order $n$ with loops; low $h$ guarantee disapproval for all such graphs; intermediate $h$ may cause approval or disapproval depending on the graph and the distribution of happiness on it. Theorem~\ref{t:Classes} determines the limits of these intervals of~$h$ for each combination of $n$ and $d$ and thereby solves Problem~\ref{pm:1} (and~\ref{pm:1'}). After proving this theorem we considered the properties of the corresponding dependencies.

Let $\hmind$ be the minimum $h$ that enables approval for given $n$ and~$d.$
At the ends of the segment $[1,n]\ni d,$ the two-stage majority becomes one-stage. Indeed, if $d=1,$ then proponents are exactly happy vertices and the final decision is made by simple majority. If $d=n,$ then the neighborhood of each vertex is the whole vertex set and simple majority on the first stage is followed by unanimity on the second one. In these cases, $\hmind=\frac{n+1}2,$ and this value assures approval for any binary labeled graph with $d\in\{1,n\}.$

For $d\approx\frac{n+1}2,$ we have $\hmind\approx\frac n4+\lD{1}$ (Proposition~\ref{p:Foot}). This means that for such $d$ and $n>9,$ there exist regular approving graphs in which only a little more than a quarter of the vertices are happy. By Proposition~\ref{p:High3} and Corollary~\ref{co:nmin}, approval with the \emph{minimum\/}
number of happy vertices $\hminn= \min_d\hmind=\lD{\frac14}(n+2+\modd{n}{4})$ is provided for such graphs by a configuration of the form $\mediastinum{S_{2h-2}}{S_{h-1}}{H_h}.$
When $d\lD{<}\frac{n+1}2,$ $\hmind$ decreases hyperbolically (with adjustment to integers) as $d$ increases (Fig.~\ref{f:ThreeDom}). Therefore, if $n$ is sufficiently large, then $\hmind$ is close to $n/4$ even for rather small~$d$ values.
\x{Say}For example, $\hmin(81,9)=23.$ This minimum is attained on a configuration $[S_{40}\;\,C_{41}]\in{\cal G}_{81|9|23}$ (Proposition~\ref{p:Low}). It is
instructive to compare this with the behavior of the ``majority of majorities" policy on $9$ disjoint voting bodies of size $9,$ which is represented by a regular graph that also belongs to ${\cal G}_{81|9},$ but requires at least $5\!\cdot\!5=25$ happy agents for\x{ the approval of a proposal} a proposal to be approved. Furthermore, the minimum $d$ such that $\hmind\le \hminn+1,$ depending on $n,$ slightly exceeds $\frac n6$ or even $\frac n8$ (Proposition~\ref{p:Foot}). When $d>\frac{n+1}2,$ $\hmind$ increases linearly with~$d$.

There are two exceptions described in Propositions~\ref{p:High3} and~\ref{p:High4} in which the fulfilment of the local and global support inequalities does not imply the possibility of approval when $d=\frac{n+1}2.$ These are the configuration classes $\calG_{5|3|2}$ and $\calG_{9|5|3}$. The regions of uniformly approving and uniformly disapproving configuration classes are symmetric with respect to the line $h = \frac n2$ on the plane with coordinates $d$ and~$h$.
Consequently, values of $d$
that enable approval with very low (about $\frac n4$) support
also require very high (about $\frac{3n}4$) support for \emph{\uniform\/}approval. The ``majority of majorities" policy with such $n$ and $d$ is $\alpha_1$-protecting and $\alpha_2$-trusting with $\alpha_1\approx\frac14$ and $\alpha_2\approx\frac34$; the corresponding security gap $2r(n,d)$ is about~$\frac n2.$

%------------------------------------------------------------
\subsection{Concluding remarks}
\label{ss:conclude}
%----------------------------------------------------

The results of this study extend the work done earlier in the literature on majority domination in graphs. They can be of interest within two real-world frameworks. Both involve two-stage decision-making and in both the second stage is a majority voting. Their difference lies in the interpretation of the first stage.

Within the first framework, an informational one, each graph vertex is an agent connected by information links with $d-1$ other agents. The agent is in one of two states, ``happy'' or ``sad,'' but when voting, it represents not only itself, but rather its neighborhood. To do this, the agent determines its vote as the majority opinion of all its neighbors (including itself). Thus, each agent aggregates local information in its vicinity\x{  environment} and such aggregated opinions\x{ states} become the votes that determine the final decision by means of the simple majority procedure.
This framework is a model of a number of settings in computer science, where the agents correspond to parallel processors, local database repositories, etc., solving, among others, error correction tasks.

The second framework represents some of the problems typical of political science.
Many socio-political decisions are made in two stages, where the first is the formation of the opinions of a number of local voting bodies, after which each such a body has one casting vote in the procedure for forming the final decision. In both stages, aggregation is carried out by a simple majority of votes.
In the simplest case, the entire set of political agents is \emph{partitioned} into a system of local voting bodies. The subcase where they are all the same size is the most democratic, since in this (and only this) subcase all voters have the same influence. This situation can be modeled by our setting involving regular graphs. Indeed, when such a graph consists of disjoint components of the same size, the vertices of any component are all proponents or all opponents, since their neighborhoods coincide. In the final voting, the number of votes of a component is equal to its size, and as their sizes are equal, this is equivalent to having one vote. Thus, we have the two-stage democratic majority mentioned above, which is $1/4$-protecting.

It is worth noting that our two-stage graph model is able to simulate a wider class of political situations in which each agent can participate in several voting bodies. In this case, the democratic requirement that the number of such bodies be the same for all agents and that they have the same size reduces to the condition that the graph determining the system of voting bodies is regular\footnote{The case where this condition is replaced by the upper boundedness of the vertex degrees is\x{ studied} considered in~\cite{LemtChe23m}.}. The main question of this study is whether\x{ plural} multiple democratic participation of each agent essentially increases the representativeness of collective decisions, and how this depends \cha{of}{on} the size of each voter's neighborhood.

In this work, it has been found that the ``curse of one-quarter'' hanging over the two-stage majority with disjoint equal-sized voting bodies is not eliminated in the case of arbitrary $d$-regular local interaction patterns, at least not for moderate~$d.$
This is the difference between our results and those of \cite{Fishburn86local,Woodall92local,Lisonek95local}, where the mentioned curse is lifted. Responsible for this is the combined effect of three differences between the two settings. Namely, in~\cite{Fishburn86local}: (1)~a special class of graphs is considered; (2)~on the second stage, only happy vertices vote, and for the approval, all of them must be proponents, instead of the condition that the proponents should make up the majority; (3)~to become a proponent, a happy vertex must have a majority of happy vertices in its neighborhood, \emph{excluding itself\xy}, thereby, its closed neighborhood must maintain a three-vote preponderance of happiness.

%----------------------------------------------------
\section*{Acknowledgements}

This paper was written while P.~Chebotarev was visiting the Weizmann Institute of Science in May--July 2022. \cha{The authors are grateful to}{We thank} Kieka Mynhardt for sharing several electronic copies of papers published in Ars Combinatoria in the 1990s, which are currently a bibliographic rarity.

The work of P.C. was supported by the European Union (ERC, GENERALIZATION, 101039692). Views and opinions expressed are, however, those of the authors only and do not necessarily reflect those of the European Union or the European Research Council Executive Agency. Neither the European Union nor the granting authority can be held responsible for them.

{\small
\bibliographystyle{elsarticle-num}%{abbrv}%{amsplain}{siam}
\bibliography{VotingRegular}

\begin{thebibliography}{10}
\expandafter\ifx\csname url\endcsname\relax
  \def\url#1{\texttt{#1}}\fi
\expandafter\ifx\csname urlprefix\endcsname\relax\def\urlprefix{URL }\fi
\expandafter\ifx\csname href\endcsname\relax
  \def\href#1#2{#2} \def\path#1{#1}\fi

\bibitem{Kalai10Noise}
G.~Kalai, Noise sensitivity and chaos in social choice theory, in: Fete of
  Combinatorics and Computer Science, Springer, 2010, pp. 173--212.

\bibitem{Nurmi99Paradoxes}
H.~Nurmi, Voting Paradoxes and How to Deal with Them, Springer, 1999.

\bibitem{Cooper13bicameral}
I.~Cooper, Bicameral or tricameral? {National parliaments and representative
  democracy in the European Union}, J. European Integration 35~(5) (2013)
  531--546.

\bibitem{Broere95majority}
I.~Broere, J.~H. Hattingh, M.~A. Henning, A.~A. McRae, Majority domination in
  graphs, Discrete Mathematics 138~(1--3) (1995) 125--135.

\bibitem{BoldiVigna02}
P.~Boldi, S.~Vigna, Fibrations of graphs, Discrete Mathematics 243~(1--3)
  (2002) 21--66.

\bibitem{GodsilNewman08}
C.~D. Godsil, M.~W. Newman, Eigenvalue bounds for independent sets, J.
  Combinatorial Theory, Series B 98~(4) (2008) 721--734.

\bibitem{HsuLin08graph}
L.-H. Hsu, C.-K. Lin, Graph Theory and Interconnection Networks, CRC Press,
  2008.

\bibitem{Melnikov98exercises}
O.~Melnikov, V.~Sarvanov, R.~Tyshkevich, V.~Yemelichev, I.~Zverovich, Exercises
  in Graph Theory, Kluwer, 1998.

\bibitem{Fishburn86local}
P.~C. Fishburn, F.~K. Hwang, H.~Lee, Do local majorities force a global
  majority?, Discrete Mathematics 61~(2--3) (1986) 165--179.

\bibitem{CockayneHedetniemi77}
E.~J. Cockayne, S.~T. Hedetniemi, Towards a theory of domination in graphs,
  Networks 7~(3) (1977) 247--261.

\bibitem{HenningHind98strict}
M.~A. Henning, H.~R. Hind, Strict majority functions on graphs, J. Graph Theory
  28~(1) (1998) 49--56.

\bibitem{KangShan00lower}
L.~Kang, E.~Shan, Lower bounds on domination function in graphs, Ars
  Combinatoria 56 (2000) 121--128.

\bibitem{Cockayne1996}
E.~J. Cockayne, C.~M. Mynhardt, On a generalisation of signed dominating
  functions of graphs, Ars Combinatoria 43 (1996) 235--245.

\bibitem{Henning96domination}
M.~A. Henning, Domination in regular graphs, Ars Combinatoria 43 (1996)
  263--271.

\bibitem{Holm01majority}
T.~S. Holm, On majority domination in graphs, Discrete Mathematics 239~(1--3)
  (2001) 1--12.

\bibitem{Ungerer96PhD}
E.~Ungerer, Aspects of signed and minus domination in graphs, Ph.D. thesis,
  Faculty of Science of the Rand Afrikaans University, Johannesburg, South
  Africa (1996).

\bibitem{HattinghUngererHenning98}
J.~H. Hattingh, E.~Ungerer, M.~A. Henning, Partial signed domination in graphs,
  Ars Combinatoria 48 (1998) 33--42.

\bibitem{Harris03PhD}
L.~M. Harris, Aspects of functional variations of domination in graphs, Ph.D.
  thesis, School of Mathematics, Statistics and Information Technology,
  University of Natal Pietermaritzburg (2003).

\bibitem{KangQiaoShanDu03}
L.~Kang, H.~Qiao, E.~Shan, D.~Du, Lower bounds on the minus domination and
  $k$-subdomination numbers, Theoretical Computer Science 296~(1) (2003)
  89--98.

\bibitem{Hattingh98majority}
J.~H. Hattingh, Majority domination and its generalizations, in: Domination in
  Graphs\/$:$ Advanced Topics, Marcel Dekker, New York, 1998, pp. 91--108.

\bibitem{ChangLiawYeh02}
G.~J. Chang, S.-C. Liaw, H.-G. Yeh, $k$-subdomination in graphs, Discrete
  Applied Mathematics 120~(1--3) (2002) 55--60.

\bibitem{ChenSong08lower}
W.~Chen, E.~Song, Lower bounds on several versions of signed domination number,
  Discrete Mathematics 308~(10) (2008) 1837--1846.

\bibitem{KangShan20signed}
L.~Kang, E.~Shan, Signed and minus dominating functions in graphs, in: Topics
  in Domination in Graphs, Springer, 2020, pp. 301--348.

\bibitem{Dunbar96minus}
J.~Dunbar, S.~Hedetniemi, M.~A. Henning, A.~A. McRae, Minus domination in
  regular graphs, Discrete Mathematics 149~(1--3) (1996) 311--312.

\bibitem{HaynesHedetniemi98book}
T.~W. Haynes, S.~Hedetniemi, P.~Slater, Fundamentals of Domination in Graphs,
  CRC Press, 1998.

\bibitem{LiuSunTian02lower}
H.-l. Liu, L.~Sun, H.-m. Tian, Lower bounds on the majority domination number
  of graphs, J. Beijing Institute of Technology 11~(4) (2002) 436--438.

\bibitem{Xing05signed}
H.-M. Xing, L.~Sun, X.-G. Chen, On signed majority total domination in graphs,
  Czechoslovak Mathematical Journal 55~(2) (2005) 341--348.

\bibitem{HarrisHattinghHenning06}
L.~Harris, J.~H. Hattingh, M.~A. Henning, Total $k$-subdominating functions on
  graphs, Australasian J. Combinatorics 35 (2006) 141–154.

\bibitem{KangShan07survey}
L.~Kang, E.~Shan, Dominating functions with integer values in graphs---a
  survey, J. Shanghai University $($English Edition$)$ 11~(5) (2007) 437--448.

\bibitem{HamidPrabhavathy16majority}
I.~S. Hamid, S.~A. Prabhavathy, Majority reinforcement number, Discrete
  Mathematics, Algorithms and Applications 8~(1) (2016) 1650014.

\bibitem{Caro18effect}
Y.~Caro, R.~Yuster, The effect of local majority on global majority in
  connected graphs, Graphs and Combinatorics 34~(6) (2018) 1469--1487.

\bibitem{Peleg02local}
D.~Peleg, Local majorities, coalitions and monopolies in graphs: {A review},
  Theoretical Computer Science 282~(2) (2002) 231--257.

\bibitem{LinialPeleg93}
N.~Linial, D.~Peleg, Y.~Rabinovich, M.~Saks, Sphere packing and local
  majorities in graphs, in: The $2$nd Israel Symp. Theory and Computing
  Systems, IEEE, 1993, pp. 141--149.

\bibitem{Peleg98size}
D.~Peleg, Size bounds for dynamic monopolies, Discrete Applied Mathematics
  86~(2--3) (1998) 263--273.

\bibitem{Auletta15minority}
V.~Auletta, I.~Caragiannis, D.~Ferraioli, C.~Galdi, G.~Persiano, Minority
  becomes majority in social networks, in: Int. Conf. on Web and Internet
  Economics, Springer, 2015, pp. 74--88.

\bibitem{Gartner18majority}
B.~G{\"a}rtner, A.~N. Zehmakan, Majority model on random regular graphs, in:
  Latin American Symp. Theoretical Informatics, Springer, 2018, pp. 572--583.

\bibitem{Avin19majority}
C.~Avin, Z.~Lotker, A.~Mizrachi, D.~Peleg, Majority vote and monopolies in
  social networks, in: Proc. 20th Int. Conf. on Distributed Computing and
  Networking, 2019, pp. 342--351.

\bibitem{Scheidler15kUnanimity}
A.~Scheidler, A.~Brutschy, E.~Ferrante, M.~Dorigo, The $k$-unanimity rule for
  self-organized decision-making in swarms of robots, IEEE Transactions on
  Cybernetics 46~(5) (2015) 1175--1188.

\bibitem{Goldberg03CollabTeleoper}
K.~Goldberg, D.~Song, A.~Levandowski, Collaborative teleoperation using
  networked spatial dynamic voting, Proceedings of the IEEE 91~(3) (2003)
  430--439.

\bibitem{Makhdoom19IoTBlockchain}
I.~Makhdoom, M.~Abolhasan, H.~Abbas, W.~Ni, Blockchain's adoption in {IoT}: The
  challenges, and a way forward, Journal of Network and Computer Applications
  125 (2019) 251--279.

\bibitem{Pourdarbani20MajorityNN}
R.~Pourdarbani, S.~Sabzi, D.~Kalantari, J.~L. Hern{\'a}ndez-Hern{\'a}ndez,
  J.~I. Arribas, A computer vision system based on majority-voting ensemble
  neural network for the automatic classification of three chickpea varieties,
  Foods 9(2), index~113, 17~pp. (2020).

\bibitem{Luccio07network}
F.~Luccio, L.~Pagli, N.~Santoro, Network decontamination in presence of local
  immunity, Int. J. Foundations of Computer Science 18~(3) (2007) 457--474.

\bibitem{Wang07majoritywatermarking}
M.-S. Wang, W.-C. Chen, A majority-voting based watermarking scheme for color
  image tamper detection and recovery, Computer Standards \& Interfaces 29~(5)
  (2007) 561--570.

\bibitem{HendrickxOlshevskyTsitsiklis11}
J.~M. Hendrickx, A.~Olshevsky, J.~N. Tsitsiklis, Distributed anonymous discrete
  function computation, IEEE Trans. Automatic Control 56~(10) (2011)
  2276--2289.

\bibitem{Cruciani21}
E.~Cruciani, E.~Natale, A.~Nusser, G.~Scornavacca, Phase transition of the
  2-choices dynamics on core--periphery networks, Distributed Computing 34~(3)
  (2021) 207--225.

\bibitem{Heilman21Stable}
S.~Heilman, Designing stable elections, Notices of the American Mathematical
  Society 68~(4) (2021) 516--527.

\bibitem{MouzonLaurent20USElect}
O.~De~Mouzon, T.~Laurent, M.~Le~Breton, D.~Lepelley, The theoretical
  {Shapley--Shubik} probability of an election inversion in a toy symmetric
  version of the {US} presidential electoral system, Social Choice and Welfare
  54~(2) (2020) 363--395.

\bibitem{KaoCouzin19modular}
A.~B. Kao, I.~D. Couzin, Modular structure within groups causes information
  loss but can improve decision accuracy, Philosophical Transactions of the
  Royal Society B 374~(1774) (2019) 20180378.

\bibitem{BottcherKernell22Condorcet}
L.~B{\"o}ttcher, G.~Kernell, Examining the limits of the condorcet jury
  theorem: {Tradeoffs} in hierarchical information aggregation systems,
  Collective Intelligence 1~(2) (2022) 26339137221133401.

\bibitem{FeixLepelleyMerlin21}
M.~Feix, D.~Lepelley, V.~Merlin, J.-L. Rouet, L.~Vidu, Majority efficient
  representation of the citizens in a federal union, in: Evaluating Voting
  Systems with Probability Models, Springer, 2021, pp. 163--187.

\bibitem{Che06ARC}
P.~Y. Chebotarev, Analytical expression of the expected values of capital at
  voting in the stochastic environment, Automation and Remote Control 67~(3)
  (2006) 480--492.

\bibitem{MaksChe20}
V.~M. Maksimov, P.~Y. Chebotarev, Voting originated social dynamics: {Quartile}
  analysis of stochastic environment peculiarities, Automation and Remote
  Control 81~(10) (2020) 1865--1883.

\bibitem{CheMal18opt}
P.~Y. Chebotarev, V.~A. Malyshev, Y.~Y. Tsodikova, A.~K. Loginov, Z.~M. Lezina,
  V.~A. Afonkin, The optimal majority threshold as a function of the variation
  coefficient of the environment, Automation and Remote Control 79~(4) (2018)
  725--736.

\bibitem{CheLog10SCW}
P.~Chebotarev, A.~Loginov, Y.~Tsodikova, Z.~Lezina, V.~Borzenko, Voting in a
  stochastic environment: {The} utility decreases as a faction expands, in: X
  Int. Meeting of the Society for Social Choice and Welfare, July 21--24, 2010,
  National Research University Higher School of Economics; {\em XI Int.
  Academic Conf. on Economic and Social Development}. Moscow, 6-8 April 2010,
  https://www.hse.ru/data/2010/04/02/1218244162/cheb.doc (in Russian), Moscow,
  2010, pp. 1--3.

\bibitem{DunbarHedeHenSla95}
J.~Dunbar, S.~Hedetniemi, M.~Henning, P.~Slater, Signed domination in graphs,
  Graph Theory, Combinatorics, and Applications 1 (1995) 311--322.

\bibitem{Lisonek95local}
P.~Lison\v{e}k, Local and global majorities revisited, Discrete Mathematics
  146~(1--3) (1995) 153--158.

\bibitem{Woodall92local}
D.~R. Woodall, Local and global proportionality, Discrete Mathematics 102~(3)
  (1992) 315--328.

\bibitem{LemtChe23m}
D.~Lemtuzhnikova, P.~Chebotarev, M.~Goubko, I.~Kudinov, N.~Shushko,
  Subdomination in graphs with upper-bounded vertex degree, Mathematics 11~(12)
  (2023) 2722.

\end{thebibliography}
}

%---------------------------
%---------------------------
\end{document}